\documentclass[preprint,3p]{elsarticle}
\biboptions{sort&compress}
\usepackage{amssymb} 
\usepackage{amsmath} 
\usepackage{amsthm} 
\usepackage{mathrsfs} 
\usepackage{bbding}
\usepackage{pifont} 
\usepackage{utfsym}

\usepackage{subcaption}
\usepackage{float}  
\usepackage{graphicx}   
\usepackage[colorlinks=true,allcolors=blue]{hyperref}   
\usepackage{tikz}   
\usetikzlibrary{shapes.geometric, arrows.meta, positioning, calc}

\tikzstyle{startstop} = [rectangle, rounded corners, minimum width=3cm, minimum height=1cm, text centered, draw=black, fill=red!30]
\tikzstyle{io} = [trapezium, trapezium left angle=70, trapezium right angle=110, minimum width=3cm, minimum height=1cm, text centered, draw=black, fill=blue!30]
\tikzstyle{process} = [rectangle, minimum width=3cm, minimum height=1cm, text centered, text width=3cm, draw=black, fill=green!30]
\tikzstyle{decision} = [rectangle, minimum width=3cm, minimum height=1cm, text centered, text width=3cm, draw=black, fill=blue!30]
\tikzstyle{arrow} = [thick,->,>=stealth]

\begin{document}
\begin{frontmatter}


\title{A hybrid PDE–ABM model for angiogenesis and tumour microenvironment with application to resistance in cancer treatment}
\author[1]{Louis Shuo Wang}
\ead{swang116@vols.utk.edu}
\author[2]{Jiguang Yu}    
\ead{zcahyuc@ucl.ac.uk}
\author[3,4]{Zonghao Liu}    
\ead{liuzonghao42@gmail.com}

\address[1]{Department of Mathematics, University of Tennessee, Knoxville, TN 37996, USA}
\address[2]{Department of Mathematics, University College London, London, WC1E  6BT, UK}
\address[3]{Clinical Oncology School of Fujian Medical University, Fujian Cancer Hospital, Fuzhou, Fujian, 350025, China}
\address[4]{Institute of Automation, Chinese Academy of Science, Beijing, 100190, China}

\begin{abstract}
The main obstacle to effective cancer treatment is the development of drug resistance, which can be divided into two categories: spontaneous and acquired drug resistance. Non-small cell lung cancer (NSCLC) is the main cause of cancer-related deaths worldwide.
A subset of lung cancer, adenocarcinomas, is characterised by mutations in the epidermal growth factor receptor (EGFR) gene.
Treatment of EGFR-mutated lung adenocarcinomas has become less effective over time due to drug resistance development, which is associated with a second mutation in the EGFR gene. An important factor in the development of cancer is angiogenesis, which is the formation of blood vessels from the existing vasculature. These newly formed blood vessels provide oxygen and nutrients to tumour cells to maintain tumour growth and proliferation. We applied a hybrid discrete-continuous (HDC) model to capture the dynamic vasculature in the tumour microenvironment (TME). In the case of pre-existing resistance, the formation of angiogenic networks creates a microenvironment that supports tumour survival and enhances drug resistance. In the case of spontaneous mutation-induced resistance, earlier and more frequent mutations confer a greater survival advantage to the tumour population. There is also a mutually reinforcing
relationship between a high proliferation rate and high resistance characteristics. These findings explain two conflicting experimental results about the second mutation in NSCLC.
\end{abstract}


\end{frontmatter}

\section{Introduction}
The main obstacle to effective cancer treatment is the development of drug resistance, which can be divided into two categories \cite{meads2009environment}: spontaneous drug resistance and acquired drug resistance. Spontaneous resistance has two forms: intrinsic spontaneous resistance, where pre-existing random genetic mutations are selected and confer resistance to the drug, and extrinsic spontaneous resistance, where the tumour microenvironment (TME) transiently protects the tumour cells from the drug. Researchers develop mathematical modelling methods to qualitatively and quantitatively study mechanisms of drug resistance and predict the eventual outcome of cancer treatment. Altrock et al. \cite{altrock2015mathematics} provides a thorough review of current mathematical models for cancers. In addition to cancer models, Yin et al. \cite{yin2019review} includes possible models for treatment resistance of solid tumours. Picco et al. \cite{picco2017integrating} integrate mathematical models for environmental-mediated drug resistance. Greene et al. \cite{greene2019mathematical} propose a model to differentiate spontaneous and induced drug resistance during cancer treatment. 

Non-small cell lung cancer (NSCLC) is the main cause of cancer-related deaths worldwide.
A subset of lung cancer, adenocarcinomas, is characterised by mutations in the epidermal growth factor receptor (EGFR) gene.
Treatment of EGFR-mutated lung adenocarcinomas has become less effective over time due to the development of drug resistance.
This resistance to drugs is associated with a second mutation in the EGFR gene, known as T790M, which is present in approximately half of patients with acquired resistance to EGFR tyrosine kinase inhibitors (TKI) \cite{chmielecki2011optimization}.

A hallmark in the development of cancer is angiogenesis, which is the formation of blood vessels from the existing vasculature, and these newly formed blood vessels penetrate the tumour, providing oxygen and nutrients to tumour cells to support further growth and proliferation of the tumour.
As tumours expand, they will exhaust oxygen and nutrients and enter the avascular stage. 
These hypoxic tumour cells in the tumour centre secrete tumour angiogenic factor (TAF) to recruit endothelial cells, thereby facilitating vessels inside the tumour. 
At this stage, tumours are vascular and enough nutrients are supplied for their expansion.
There is an extensive body of work on the mathematical modelling of angiogenesis and its impact on tumour growth \cite{billy2009pharmacologically,jackson2000mathematical,bodzioch2021angiogenesis,sun2016mathematical,spill2015mesoscopic}. Balding and McElwain \cite{balding1985mathematical}, Byrne and Chaplain \cite{byrne1995mathematical}, and Pillay et al. \cite{pillay2017modeling} each formulalted reaction diffusion systems for tumour-induced vessel formation and numerically validated the travelling wave solutions for tip and vessel dynamics. Anderson and Chaplain \cite{anderson1998continuous} pioneered a hybrid discrete-continuous approach for tumour angiogenesis and set the rules for endothelial cell movement and proliferation, and vessel branching and anastomosis. Their model numerically reproduced the brush-border effects of the vasculature network. Flandoli et al. \cite{flandoli2023mathematical} employed stochastic Poisson processes for cell cycles and vessel dynamics, and introduced the Lenard-Jones potential for the repulsion-attraction interaction between tumour cells.

Agent-based modelling approaches have become a powerful tool in computational biology studies, for their capabilities of simulating the diversity and complexity of the interactions within each cell, different cells within a cell population, or across populations, and between cells and their local microenvironment. They are well suited to studying biological systems with high levels of heterogeneity, which is difficult to achieve with other methods. Soheilypour and Mofrad \cite{soheilypour2018agent} summarised the capabilities of agent-based models in molecular system biology modeling and exploring the dynamics of molecular systems and pathways in health and disease. Wang et al. \cite{wang2015simulating} reviewed the scale ranges and hybrid modeling aspects of agent-based models and their applications in molecular signalling, cellular metabolism, phenotypic changes, angiogenesis, microenvironment, metastasis, cancer stem cells, cancer treatment and immune response. Anderson \cite{anderson2007hybrid} used hybrid discrete-continuous techniques to study the impact of the tumor microenvironment on solid tumour growth and invasion, and they found that harsher microenvironments select for more invasive and aggressive tumour phenotypes. Gevertz et al. \cite{gevertz2015emergence} employed the spatial agent-based model to explore the importance of microenvironmental niche on the emergence of anti-cancer drug resistance. They found the complex interaction between the niches and the response of tumour cells to the drug on the development of drug resistance. Picco \cite{picco2024role} investigated the role of environmentally mediated drug resistance in facilitating the spatial distribution of residual disease. Yang et al. \cite{yang2023multiscale} used multiscale modelling techniques of drug resistance and their findings underscore the dual roles of intrinsic genetic mutations and extrinsic microenvironmental adaptations in steering tumor growth and drug resistance.

Thus, this research is conducted to study the role of angiogenesis in the development of drug resistance, where pre-existing drug resistance and induced drug resistance are considered.
We propose a hybrid discrete-continuous model (HDC) to capture the dynamic vasculature in the tumour microenvironment (TME),
where agent-based models characterise the tumour and endothelial cells, and reaction-diffusion equations describe the drug and oxygen fields. 
The remainder of the paper is structured as follows. Section~\ref{sec:modeling} presents the mathematical modeling of the hybrid discrete-continuous model,
including the PDE model for TME and the agent-based model for tumour cell and angiogenesis dynamics. 
The numerical results are presented in Section~\ref{sec:numerical},
where we investigate the effects of angiogenesis on drug resistance and the optimal drug delivery strategies to overcome drug resistance.
Section~\ref{sec:biology} discusses the biological implications of our work, where our result successfully explains two contradictory experiments in the literature.
We list future directions and potential extensions of our work in Section~\ref{sec:future}.
Finally, Section~\ref{sec:conclusion} concludes the paper with a summary of the findings and future work.

\section{Modeling}\label{sec:modeling}
\subsection{PDE Model for Tumour Microenvironment}

We aim to study, under tumour-induced angiogenesis, the spontaneous and induced mechanisms of drug resistance of cancer, and we hope to use the model to find the optimal drug delivery strategies to overcome drug resistance. 

We adopt the modelling in \cite{anderson1998continuous} 
to capture the dynamic vasculature in the TME. There are two contributions 
to the flux of endothelial cells $n$: random diffusion, chemotaxis, and haptotaxis: \begin{align*}
J_{n}=J_{\mbox{diff}}+J_{\mbox{chemo}}
 \end{align*} 
The random diffusion is given by:
\begin{align*}
J_{\mbox{diff}} = -D_{n} \nabla n
\end{align*}
where $D_{n}$ is the diffusion coefficient of endothelial cells. The chemotaxis motion is stimulated by tumour angiogenic factor (TAF) $c$
and the chemotaxis flux is taken to be $J_{\mbox{chemo}}=\chi(c)n\nabla c $, where $\chi (c)$ is the chemotaxis function taken as 
$\chi (c)=\dfrac{\chi _{0}k_{1}}{k_{1}+c} $ with $\chi _{0}>0$ being the chemotaxis coefficient and $k_{1}>0$ being a positive constant. 
The equation for endothelial cells is given by:
\begin{align}\label{eq1.1}
\dfrac{\partial n}{\partial t}=D_{n}\mathop{{}\bigtriangleup}\nolimits n-\nabla \cdot \left( \chi (c)n \nabla c \right).   
 \end{align} 

The TAF concentration $c$ is assumed to satisfy the following equation:
\begin{align}\label{eq1.2}
\dfrac{\partial c}{\partial t}=D_{c}\mathop{{}\bigtriangleup}\nolimits c - \xi _{c}c + \eta \sum\limits_{a \in \Lambda ^{h}(t)} \chi _{a}-\lambda c\sum\limits_{v\in V_{t}} \chi _{v}, 
 \end{align}  
where $D_{c}$ is the diffusion rate of TAF, $\eta$ is the production rate of TAF by tumour cells 
with $\chi _{a}$ being the Dirac function concentrated at a hypoxic tumour cell site $a \in \Lambda^{h}_{t}$, 
and $\lambda $ is the uptake rate by endothelial cells with $\chi _{v}$ being the Dirac function concentrated at the vessel sites $v \in V_{t}$. 

The vasculature or blood vessels can have an essential impact on drug delivery through blood vessels, thus facilitating the modelling of the heterogeneity
of drug distribution in the TME. We use $d$ to denote the drug concentration, and its delivery and cellular uptake are assumed to satisfy the following equation:
\begin{align}\label{eq1.3}
\dfrac{\partial d}{\partial t}(x,t)=D_{d}\mathop{{}\bigtriangleup}\nolimits d(x,t)-\xi_{d}d(x,t)-\rho _{d}d(x,t)\sum\limits_{a \in \Lambda _{t}} \chi _{a}(x,t)+S_{d}(t)\sum\limits_{v\in V_{t}} \chi _{v}(x,t),  
\end{align}  
where $D_{d}$ is the diffusion rate of the drug, $\xi_{d} $ is the decay rate, 
$\rho _{d}$ is the cellular uptake rate by tumour cells, 
and $S_{d}(t)$ is the drug supply rate.

The definition of the characteristic functions $\chi _{a}, \chi _{v}$ depends on the cell radius:
\begin{align*}
\chi _{a}(x,t)=\begin{cases}
1 &  \mbox{ if } \left\lVert x-a^{\left( X,Y \right) }(t) \right\rVert \leq R_{c} \\ 
0 & \mbox{ otherwise} 
\end{cases} \quad \mbox{ and } \quad \chi _{v}(x,t)=\begin{cases}
    1 &  \mbox{ if } \left\lVert x-v^{\left( X,Y \right) }(t) \right\rVert \leq R_{c} \\ 
    0 & \mbox{ otherwise} 
    \end{cases}
 \end{align*} 

In addition to drugs, blood vessels provide nutrients, such as oxygen, for tumour cells. In the avascular phase of
tumour invasion, the rapid proliferation of tumour cells will exhaust nearby oxygen in the TME, thus rendering tumour cells hypoxic. 
Hypoxic cells in the necrotic core secrete TAF, which diffuses in TME and attracts the chemotaxis movement of endothelial cells toward tumour cells. 
Once endothelial cells reach the tumour, new blood vessels within the tumour begin to grow, and the vascular phase of
tumour invasion starts. These new vasculatures within the tumour provide oxygen to the tumour and promote further
proliferation and metastasis of tumour cells. We assume the oxygen concentration $o$ changes according to the following rule:
\begin{align}\label{eq1.4}
\begin{split}
\dfrac{\partial o}{\partial t}(x,t)&=D_{o}\mathop{{}\bigtriangleup}\nolimits o(x,t)-\xi_{o}o(x,t)
-\rho _{o} \sum\limits_{a \in \Lambda _{t}} \chi _{a}(x,t)  \\ 
                                   &+S_{o}\left( o_{\max}-o \right) \sum\limits_{v \in V_{t}}\chi _{v}(x,t), 
\end{split}
\end{align} 
where $D_{o}$ is the diffusion rate of the oxygen, $\xi_{o}$ is the decay rate, 
$\rho _{o}$ is the cellular uptake rate by tumour and endothelial cells, and $S_{o}$ is the supply rate at blood vessel sites, and $o_{max}$ is the maximum oxygen concentration in TME with the factor $o_{\max}-o$ incorporating the fact 
that less oxygen is released through the vessel walls at high environmental oxygen concentrations.

All PDEs in this section are updated using the finite-difference method with the forward Euler scheme.
We necessarily impose homogeneous Neumann boundary conditions to model the closed no-flux system:
\begin{align*}
\vec{n}\cdot \left( -D_{n}\nabla n+\chi (c)n \nabla c\right) =0,\\
\vec{n}\cdot \left( D_{c}\nabla c\right) =0,\\
\vec{n}\cdot \left( D_{d}\nabla d\right) =0,\\
 \end{align*} 
on the boundary of the spatial domain with $\vec{n}$ being the outward unit normal vector.

\subsection{Agent-based Model for Tumour Cell Dynamics}

To describe the dynamics of tumour and angiogenesis, we first introduce some notation. 
Let $\Lambda _{t}$ be the collection of tumour cells at time $t$, and $\Lambda _{t}^{n}, \Lambda _{t}^{h}$ be the collection of normoxic and hypoxic tumour cells at time $t$, respectively.
The total number of tumour cells at time $t$ is $N_{t}=\left\lvert \Lambda _{t} \right\rvert $, and the total number of normoxic and hypoxic tumour cells at time $t$ are $N_{t}^{n}=\left\lvert \Lambda _{t}^{n} \right\rvert $ and $N_{t}^{h}=\left\lvert \Lambda _{t}^{h} \right\rvert $, respectively.

Here comes the agent-based model for tumour cell dynamics. Each tumour cell $a \in \Lambda _{t}$ has a unique index $a=(k,i_{1},\ldots ,i_{n})$ with $1\leq k\leq N_{0}$ and $i_{1},\ldots i_{n}\in \left\{ 1,2 \right\} $. Generally speaking, $n$ represents the generation after the initial $N_{0}$ cells to which this tumour cell $a$ belongs. The initial
$N_{0}$ cells belong to the generation $0$ ($n=0$) and have a simple index $\mathrm{id}_a=\left( k \right) $. After each division of an $n$-th generation cell $\mathrm{id}_a=\left( k,i_{1},\ldots,i_{n} \right) $,  
the indices of the two daughter cells are inherited from the mother cell:
\begin{align*}
\mathrm{id}_{a\prime} =\left( k,i_{1},\ldots ,i_{n},i_{n+1} \right),\,\, i_{n+1}=1,2. 
\end{align*}
We regard each tumour cell as an agent associated with a unique set of regulated properties: 
\begin{align}\label{eq1.6}
a=\left\{ a^{\left( X,Y \right) }(t), a^{o}(t), a^{d}(t), a^{dam}(t), a^{death}(t), a^{age}(t), a^{mat}, \mathrm{id}_a \right\}. 
 \end{align} 
$a^{\left( X, Y \right) }$ defines the position of the tumour cell, and the following rule updates it:
\begin{align*}
a^{\left( X,Y \right) }(t+\Delta t)=a^{\left( X,Y \right) }(t)+\left( \Delta X, \Delta Y \right). 
\end{align*}  
Each cell can inspect its local neighbourhood and sense extracellular concentrations of oxygen $o$ and drug $d$.
The level of oxygen taken up and used by tumour cells is:
\begin{align*}
a^{o}(t+\Delta t)= o \left( a^{\left( X,Y \right) }(t),t \right) ,
 \end{align*} 
There are two critical oxygen concentration thresholds $o_{apop}<o_{hyp}$ such that 
tumour cells with cellular oxygen concentration $a^{o}>o_{hyp}$ are normoxic cells and $a^{o}\leq o_{hyp}$ hypoxic cells.
Furthermore, tumour cells with $a^{o}\leq o_{apop}$ are insufficient to maintain survival and undergo apoptosis, 
and we remove these cells from our square immediately. We differentiate tumour cells of normoxic and hypoxic types
as $a \in \Lambda _{t}^{n}$ and $a \in \Lambda _{t}^{h}$. 
  
The amount of drug taken up by cells $a^{d}$ is determined by:
\begin{align*}
a^{d}(t+\Delta t)=a^{d}(t)+ d \left( a^{\left( X,Y \right) }(t),t \right) \Delta t,
 \end{align*} 
We update the drug-induced cellular DNA damage $a^{dam}$ as follows:
\begin{align*}
a^{dam}(t+\Delta t)=a^{dam}(t)+ d \left( a^{\left( X,Y \right) }(t),t \right) \Delta t-p_{r}a^{dam}(t)\Delta t,
 \end{align*} 
where $p_{r}$ is the cellular repair rate of the DNA damage. 
The death threshold $a^{death}$ is the maximum tolerated level of DNA damage,
and the cell will die if the level of DNA damage $a^{dam}$ exceeds this threshold. In the pre-existing mechanism case, the death threshold is fixed for all cells and is set to be higher for resistant cells. 
We assume, for simplicity, that $a_{R}^{death}=Th_{multi}\times a_{S}^{death}$ with $Th_{multi}$ being a positive integer. 
The cell's current age will increase if tumour cells remain in the normoxic state, and we update the current cell age as:
\begin{align*}
a^{age}(t+\Delta t)=\begin{cases}
a^{age}(t)+\Delta t & \mbox{ if } a^{o}(t)>o_{hyp} \\
a^{age}(t)          & \mbox{ otherwise }
\end{cases}.
 \end{align*} 
The cell maturation time $a^{mat}$ is dependent upon normoxic cell division rate $\alpha _{n}$: $a^{mat}=\dfrac{\log(2)}{\alpha _{n}} $.  
When the cell age of normoxic cells reaches the cell maturation time $a^{mat}$, normoxic cells are ready to
proliferate and divide into two daughter cells if the neighbourhood of the mother cell has at least one empty space, and proliferation will be suppressed if all neighbourhood squares are occupied by tumour cells.
Upon division of $a$, two daughter cells are created simultaneously, with one daughter cell 
occupying the mother cell position and the other daughter cell occupying a neighboring square; that is:
\begin{align*}
&a_{1}^{\left( X,Y \right) }(t)=a^{\left( X,Y \right) }(t), \\
&a_{2}^{\left( X,Y \right) }(t)=a^{\left( X,Y \right) }(t)+\left( \Delta X_{k}, \Delta Y_{k} \right),
\end{align*} 
where $(\Delta X_{k},\Delta Y_{k})$ is the random movement into the neighboring position of the mother cell position.
The level of oxygen content in the daughter cells is determined by the oxygen concentrations in their locations.
The current cell age of the daughter cells $a_{1}^{age}(t),a_{2}^{age}(t)$ after division is set to zero. 
The level of cellular DNA damage and the accumulated drug of a mother cell are symmetrically segregated between the two daughter cells:
$a_{1}^{dam}(t)=a_{2}^{dam}(t)=\dfrac{1}{2}a^{dam}(t)$, $a_{1}^{d}(t)=a_{2}^{d}(t)=\dfrac{1}{2}a^{d}(t)$.
Other properties of daughter cells are determined by mutation algorithms, which will be discussed in a later section.
The initial properties of the $0$-th generation tumour cell are set to be: $(k)(0)=\left\{ \left( X_{k},Y_{k} \right), 
o \left( a^{\left( X, Y \right) }(0),0 \right), 0, 0, T_{k}, N_{k}, M_{k}, (k)  \right\} $,
where $M_{k} = \wp _{age}$ with $\wp _{age}$ being the average maturation age, which follows from a
uniform distribution $U[0.9, 1.1] \mbox{ days}$, and $N_{k}$ follows from a uniform distribution $U[0, M_{k}]$.  
$T_{k}=Th_{death}$ for all cells in the acquired mechanism case and all sensitive cells in the pre-existing mechanism case. 
$T_{k}=Th_{multi}\times Th_{death}$ for all resistant cells in the case of the pre-existing mechanism. 

\subsection{Local density of tumour cells}

We introduce the local density of the tumour cells to model the crowding effect of the tumour cells. 
We consider the indicator function of a ball of raius $R_{c}$ centered at $0$, $B_{R_{c}}$:
\begin{align*}
W(x)=\chi _{B_{R_{c}}}(x),
 \end{align*}    
and the function $F$ to indicate the number of tumour and endothelial cells around $x$:
\begin{align*}
F(x,t)=\sum\limits_{\tilde{a} \in \Lambda _{t}} W \left(x-\tilde{a}^{\left( X,Y \right) }(t) \right)+\sum\limits_{v \in V_{t}} W \left(x-v^{\left( X,Y \right) }(t) \right).
 \end{align*}  
Let $F_{max}$ be the maximum number of tumour and endothelial cells that can be accommodated in the vicinity of a tumour cell $a$,
then the tumour cell $a$ ceases to proliferate under the crowding effect, i.e., $F \left( a^{\left( X,Y \right) }(t),t \right)\geq F_{\max} $. 

\subsection{Mutation}
Some works incorporate different mutation mechanisms. To mention a few, the random mutation algorithm presented
in \cite{anderson2007hybrid} employs a selection process involving 100 distinct phenotypes, 
each assigned an equal probability of selection during mutation events. Furthermore, the study examines a linear mutation algorithm in which,
after mutation, cells acquire the phenotype at the next level, characterised by increased resistance and aggressiveness. 
Although the linear mutation algorithm is regarded as more biologically plausible than the random mutation algorithm, due to the latter's potential to induce abrupt changes in trait values, it remains insufficiently aligned with biological reality. 
The inadequacy arises from its fixed mutation direction, which consistently promotes progress towards more resistant and aggressive phenotypes. The linear algorithm is predicted to always result in the most aggressive phenotypes taking over the entire population in sequence, while evolutionary selection pressures imposed by the microenvironment will not play a role.

Therefore, to rectify these limitations, we propose our consecutive random mutation algorithm designed to mitigate sudden trait changes
while allowing for a nondirectional approach to mutation. 
The term 'consecutiveness' in our mutation algorithm refers to the gradual alteration of each trait's value through multiplication by a random number, contrasting with other random algorithms that may induce abrupt changes.
When tumour cells divide, daughter cells may acquire mutations. We assume a constant mutation rate and
propose a consecutive random mutation algorithm in which, each time a mutation occurs, three random numbers are generated from
a uniform distribution $U[0.7,1.7]$. The value of each trait is then multiplied by the corresponding random number, with the restriction that the new value must fall within the range $[0.5x,4x]$ where $x$ is the trait value of non-mutated cells.

\subsection{Equations of Tumour Cell Mechanics}
Each tumour cell in our model is an individual agent represented by the coordinates of its nucleus $a^{\left(X, Y\right)}(t)$ and a fixed cell radius $R_{C}$.
The position of a tumour cell $a$ is subject to the following equation:
\begin{align*}
da^{\left( X,Y \right) } = \varepsilon _{1}dW_{t}^{a},
\end{align*} 
where $\left( W_{t}^{a} \right)_{a\in \Lambda_{t}} $ is a sequence of independent Brownian motions, 
and $\varepsilon _{1}$ is the noise intensity.
We update the position by the Euler-Maruyama scheme:
\begin{align*}
a^{\left( X,Y \right) }(t+\Delta t)=a^{\left( X,Y \right) }(t)+ \sqrt{\Delta t} \varepsilon _{1} Z_{t}^{a},
 \end{align*} 
where $\left( Z_{t}^{a} \right)_{a\in \Lambda_{t}} $, for each fixed $t\geq 0$, are 
a sequence of i.i.d. random variables drawn from the standard normal distribution.

When a dividing cell gives rise to two daughter cells, the repulsive-attractive forces between the daughter cells are activated since they are placed at a distance smaller than the cell diameter. In addition, this may result in the placement of daughter cells near other tumour cells. Therefore, multiple repulsive and attractive forces will be applied until the whole tumour cluster reaches an equilibrium configuration under cell interactions.

\subsection{Modeling of Angiogenic Network}
We base our model on the assumption that the motion of an individual endothelial cell located at the tip of a capillary sprout
governs the motion of the entire sprout, since the remaining endothelial cells that line the wall of the sprout are contiguous.
In our schematic representation of blood vessels, the walls of the blood vessels consist of a layer of contiguous endothelial cells, and thus, each side of the vessel wall is represented
by one line of endothelial cells. Each endothelial cell is a long, thin cell, and the length of one endothelial cell
and the width of the blood vessel are set to be the diameter of the cell $2R_{c}$. The endothelial cell at the sprout tip has the shape of a semicircle arc with radius $R_{c}$, and each endpoint of the arc is connected to one line of endothelial cells.

For modelling purposes, we discretised the spatial domain into small squares, with the side length of each square being the diameter of
the endothelial cell $2R_{c}$. We model the endothelial cell at the sprout tip as an agent $b\in T_{t}$ occupying one square:
\begin{align*}
b=\left\{ b^{\left( X,Y \right) }(t), b^{age}(t), \mathrm{id}_b \right\},
 \end{align*} 
 where $b^{\left( X,Y \right) }$ is the center position of the square occupied by the tip cell, and it is introduced to 
 keep track of the position of the tip cell, $b^{age}$ is the age of the tip cell counted from the moment when it is branched from an existing tip cell, 
$\mathrm{id}_b$ as a tip cell property is the index of the tip cell $b$ that is updated in the same way as the tumour cell index, and $T_{t}$ is the collection of all tip cells at time $t$.  
The trajectory of the tip cell is given in detail in \ref{app:a1}, and we model the angiogenesis network $A_{t}$ at the time $t$
by the trajectories of all tip cells up to time $t$:
\begin{align*}
A_{t}=\bigcup\limits_{b\in T_{t}} \left\{ b^{\left( X,Y \right) }(s): 0\leq s\leq t \right\}. 
 \end{align*}  
 We update the age of the tip cell $b$ as:
 \begin{align*}
    b^{age}(t+\Delta t)=b^{age}(t)+\Delta t.
  \end{align*}  
We split the blood vessels into discrete agents and approximate vessels in the following way: Each square having a non-empty intersection with the angiogenic network $A_{t}$ is considered an agent $v=\left\{ v^{\left( X,Y \right) }(t) \right\} \in V_{t}$ and we call them vessel agents. $v^{\left( X,Y \right) }$ is the centre position of the square occupied by the vessel agent $v$, and $V_{t}$ is the collection of all the vessel agents at time $t$. In the vascular phase of tumour invasion, blood vessels grow inside the tumour and occupy space. The vessel agents have repulsive-attractive interactions
with tumour cells, and we assume that each vessel agent $v\in V_{t}$ is a circle with radius $R_{c}$ and centre position $v^{\left( X,Y \right) }(t)$. 

\subsection{Rules for Branching, Anastomosis, Regression, and Proliferation}
During angiogenesis, there are two processes: the generation of new sprouts (sprout branching) and the connection of two existing sprouts (sprout anastomosis).

We will assume that the generation of new sprouts (branching) occurs only from existing sprout tips and that the newly formed sprouts are unlikely to branch immediately.
From these assumptions, we obtain the following necessary conditions for branching.
\begin{enumerate}[(i)]
    \item The age of the current sprout $b^{age}(t)$ is greater than some threshold branching age $\psi $, i.e., new sprouts must mature for a length of time
at least equal to $\psi $ before being able to branch.  

    \item There is sufficient space locally for a new sprout to form; that is, branching into a space occupied by another sprout is not allowed.
\end{enumerate} 
When the above conditions are satisfied, we assume the branching event of each tip cell $b\in T_{t}$ follows from a Poisson distribution with intensity:
\begin{align*}
    \lambda_{br}=c_{br}\, c \left( b^{\left( X,Y \right) }(t),t \right) / \left\lVert c(\cdot,t) \right\rVert _{\infty } ,
\end{align*} 
where $c_{br}$ is a positive constant, and $\left\lVert c(\cdot,t) \right\rVert _{\infty } $ is the maximum value of the TAF field at the current time step. 
When a tip cell $b$ branches, two new tip cells $b_{1},b_{2}$ are created, with one tip cell occupying the position of the old tip cell and the other tip cell occupying an orthogonal neighbouring square. The age of $b_{1},b_{2}$ is set to zero.

We can capture the anastomosis process with our discrete model. As the sprout progresses towards the tumour, the endothelial cell
at the sprout tips can move to any of the four orthogonal neighbours on the discrete grid. If in one of these moves another sprout is encountered, then anastomosis can occur. As a result, only one of the original sprouts continues to grow, the choice of which is purely random.

We modelled the endothelial cell proliferation process, assuming that some of the cells behind the tip of the sprout divide into two daughter cells every 18 hours.
This has the effect of increasing the length of a sprout by approximately one cell length every 18 hours. 

\section{Numerical Results}
\label{sec:numerical}
\subsection{Initial TAF Field}
\begin{figure}[htbp]
    \centering
    \hspace*{-2cm}\includegraphics[width=1.25\linewidth]{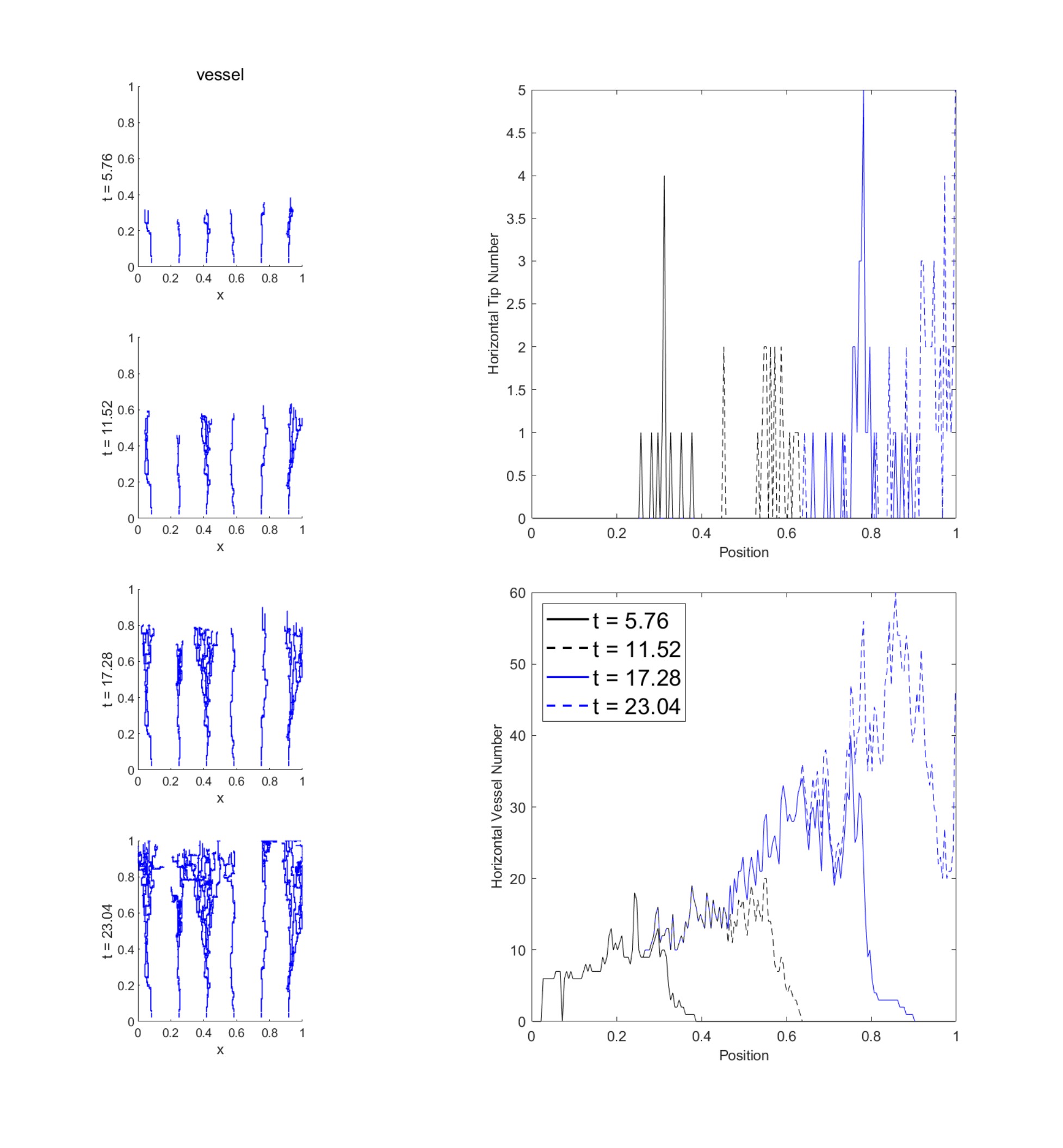}
    \caption{Vascular network up to $t=23.04$}
    \label{fig:notumour}
\end{figure}

We first set the initial TAF concentration to be a linear function in the $y$ direction: $c(x,y,0) \sim y $, and this initial setting
is motivated by the work \cite{pillay2017modeling}. The timescale for TAF diffusion is much faster than the timescale for
vascular network formation, and this justifies the application of a quasi-steady state profile for initial TAF concentration.
The TAF field is dominated by
\begin{align*}
& \mathop{{}\bigtriangleup}\nolimits c = 0, \mbox{ for } x,y\in [0,1], \\ 
& c(x,0) = 0, c(x,1) = k, \mbox{ for } x\in [0,1], \\ 
& \dfrac{\partial c}{\partial x}(0,y) = \dfrac{\partial c}{\partial x}(1,y) = 0, \mbox{ for } y\in [0,1],
 \end{align*} 
where $k$ is some positive constant to be specified later. This mimics the situation where a line source of tumour cells is located
at $y=1$ and maintains a constant concentration of TAF \cite{anderson1998continuous}. The solution to the above TAF problem is \begin{align*}
c(x,y) = ky.
 \end{align*} 
We must choose the value of $k$ so that our vascular growth is in agreement with some existing experimental data. 
In \cite{balding1985mathematical}, $14$ days is a typical value in many experiments for the time required to complete
vascularization when the tumour is implanted at a distance of $2$ mm from the limbus. The first six tip cells are distributed at $(x,y) = \left( \dfrac{2i-1}{12},0.05 \right) $. After trial and error, we find that the appropriate value for $k$ is $k = 5$, under which the time to complete vascularization in our simulation is $14.08$ days, and the final pattern of vessels in Figure~\ref{fig:notumour} is
qualitatively similar to the result in Fig.~$11$ of \cite{anderson1998continuous}, where in our final vasculature we observed the effects of "brush board": a short distance from the tumour, the frequency of branching increases dramatically.

We count the horizontal vessel number at each constant $y$-slice of our 2D domain, and in Figure~\ref{fig:notumour} we observe the travelling
wave phenomenon for tip cells and endothelial cells in the propagation of the vessel front predicted by many other works \cite{pillay2017modeling,balding1985mathematical,byrne1995mathematical}. 

Extensive simulations (results not shown) indicate that almost all anastomosis events occur when a tip cell encounters the endothelial cells
in its sprout (self-loops), creating a stunted sprout. Many tip cells are annihilated through these self-loops, and this
accounts for the deviations between our vessel network and many well-formed vessel networks in the literature. 
Such self-loops are inefficient as they will not perfuse or connect the tumour to the limbus (blood supply). 
In addition, in practice, it is unlikely that many stunted sprouts occur and, if they do, most will likely regress \cite{carmeliet2011molecular}.

In our code, we consider a sprout to belong to some tip cell, say $b_{i}$, if this sprout is generated by the trajectory of $b_{i}$ 
during the time period $[t_{1}, t_{2}]$, where $t_{1}$ is the time when $b_{i}$ is generated by branching or $t_{1}=0$ if $b_{1}$ is the initial tip cells in our simulation,
and $t_{2}$ is the time when $b_{i}$ branches to give rise to two new tip cells if no anastomosis occurs before branching, or $t_{2}$ is the time when
$b_{i}$ anastomoses with other sprouts or tip cells.

\subsection{Resistance Mechanisms}
\label{sub4.2}
There are two treatment protocols in our simulation: continuous low-dose infusion and periodically pulsed high-dose infusion with on- and off-treatment. The two strategies are compared under the assumption that the total dose in a treatment period is the same.
The continuous infusion is given by:
\begin{align*}
S_{d}(t) = d_{c}, 
 \end{align*} 
 while the pulsed infusion is given by:
 \begin{align*}
 S_{d}(t) = \begin{cases}
 d_{p} & t\in t_{init}+i(t_{on}+t_{off})+[0, t_{on}] \\ 
 0 & t\in t_{init}+i(t_{on}+t_{off})+(t_{on}, t_{on}+t_{off}]
 \end{cases} 
  \end{align*} 
for $i\in \mathbb{Z}_{\geq 0}$, where $t_{init}$ is the treatment starting time, and $t_{on},t_{off}$ are the time duration of treatment and treatment holiday, respectively.
The total dose in a treatment period is $d_{c}(t_{on}+t_{off})$ and $d_{p}t_{on}$ in these two cases, and we deduce that 
\begin{align*}
d_{c} = d_{p}\dfrac{t_{on}}{t_{on}+t_{off}}
\end{align*} given the above assumption.  
We incorporate three types of resistance mechanisms into our model: pre-existing resistance, spontaneous mutation, and drug-induced mutation.
Initially, there are $1 \%$ resistance cells whose death threshold is $Th_{multi}>1$ times that of sensitive cells in the pre-existing setting.
We consider a random mutation algorithm. Each time a mutation occurs, three random numbers will be generated from the uniform distribution $U[0.7,1.7]$
and the value of each trait (death threshold, oxygen consumption rate, and proliferation rate) will be multiplied by the corresponding random number
under the constraint that the new value falls within the range $[0.5x,4x]$ where $x$ is the trait value of non-mutated cells. 
The effect of drug-induced resistance is modelled by the drug-dependent mutation algorithm, where, under the influence of the drug, this mutation rate is proportional to the drug concentration.

\subsection{Vascularization}
\begin{figure}
    \centering
    \hspace*{-2cm}\includegraphics[width=1.25\linewidth]{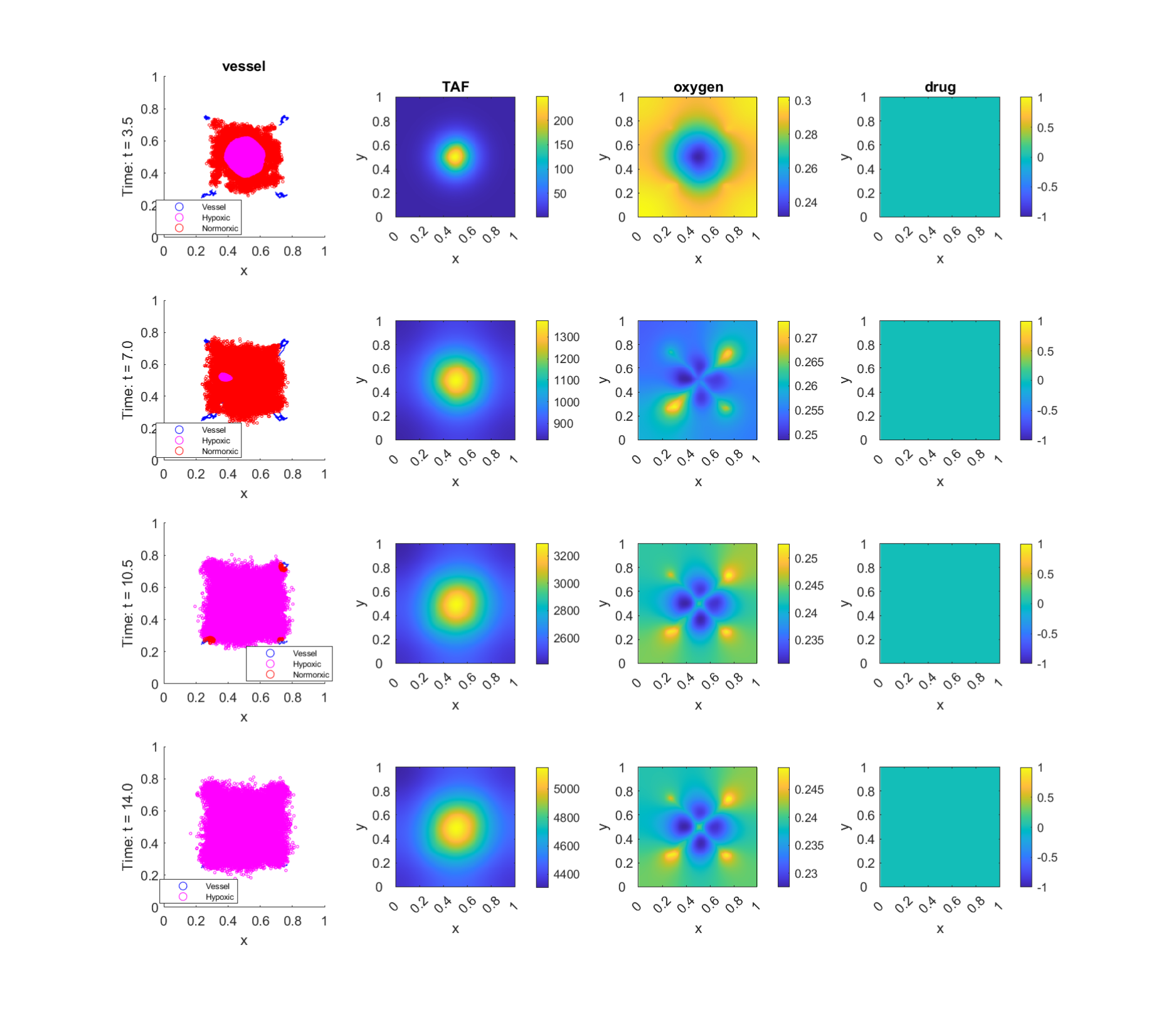}
    \caption{Vascularization}
    \label{fig:avascular1}
\end{figure}

\begin{figure}
    \centering
    \hspace*{-2cm}\includegraphics[width=1.25\linewidth]{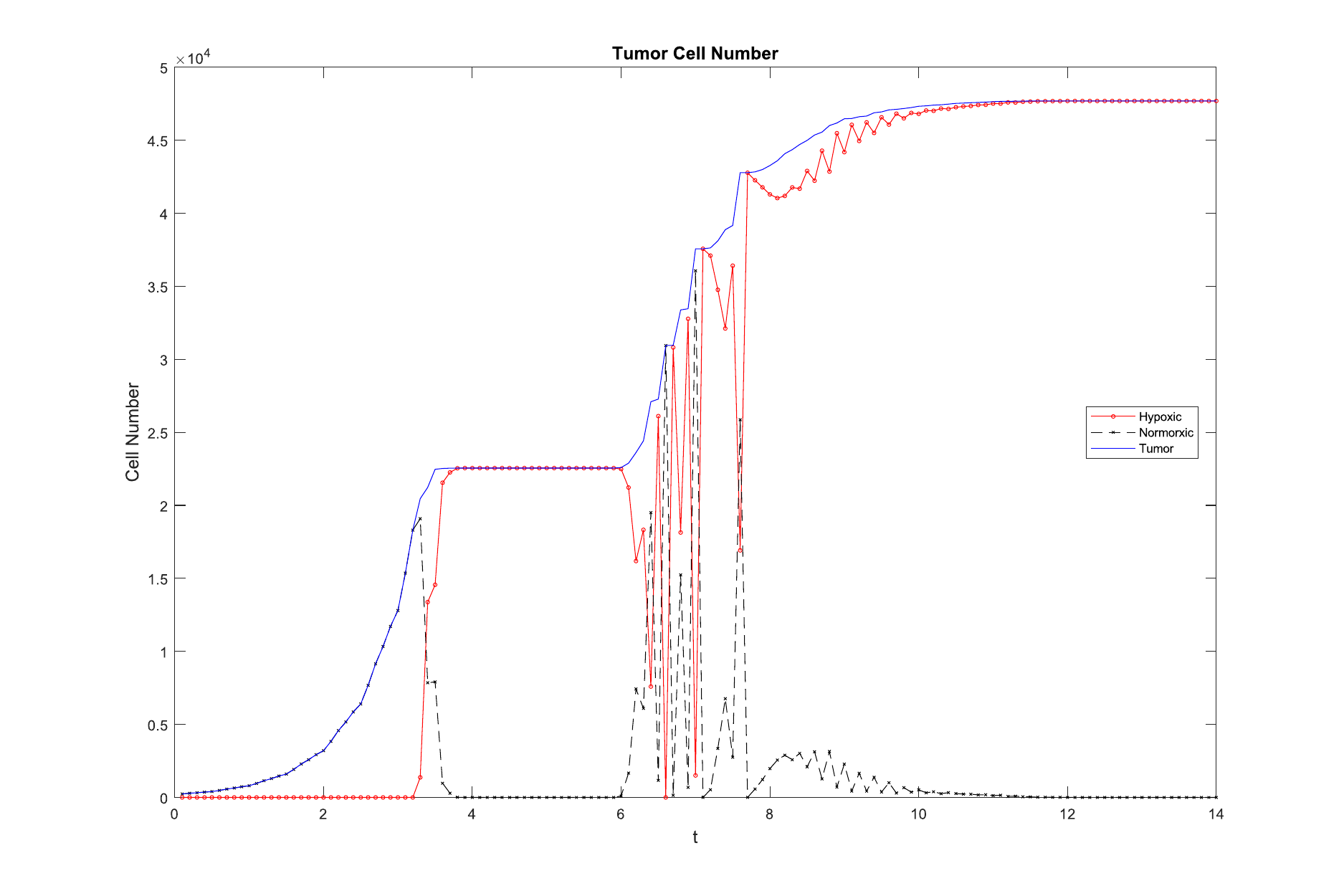}
    \caption{Tumor number under vascularization}
    \label{fig:avascular2}
\end{figure}

As depicted in Figure~\ref{fig:avascular1} and Figure~\ref{fig:avascular2}, the tumour population first experiences exponential growth with a marked expansion of the tumour mass until it exhausts the available oxygen, when the hypoxic region is formed inside the tumour and tumour proliferation stops.
After this avascular phase, hypoxic cells secrete TAF that diffuses into surrounding tissues, recruits endothelial cells,
and stimulates the formation of new blood vessels in the tumour mass. 
The TAF gradient is steep in the centre of the tumour and shallow elsewhere.
These blood vessels penetrate the tumour and provide the necessary oxygen, and this is the vascular stage of tumour development.
At the beginning of this stage, there is a strong oscillation between normoxic and hypoxic cells, and this oscillation is due to the assumption that the oxygen consumption rate of hypoxic cells is slower compared to that of normoxic cells. Then, when normoxic cells dominate, they consume oxygen faster, leading to hypoxia; while hypoxic cells dominate, they consume less oxygen, resulting in normoxia. The oscillation continues until the tumour population reaches a new carrying capacity, and in this steady state, with hypoxic cells predominating the tumour population, the sum of oxygen consumption and decay rates equals the supply rate, and the oxygen distribution remains static, and this is the mature stage of tumour development. 

We used this vascularised tumour model to study the effect of environmental factors and resistance mechanisms on tumour growth and response to treatment.

\subsection{No Resistance}
The damage clearance rate $p_{r}=0.2$. If resistance does not occur and there are no resistant cells before treatment, the tumour will be eliminated after a few cycles of treatment, where constant infusion of the drug is applied from $t = 14$. 
\begin{figure}[htbp]
    \centering
    \hspace*{-2cm}\includegraphics[width=1.25\linewidth]{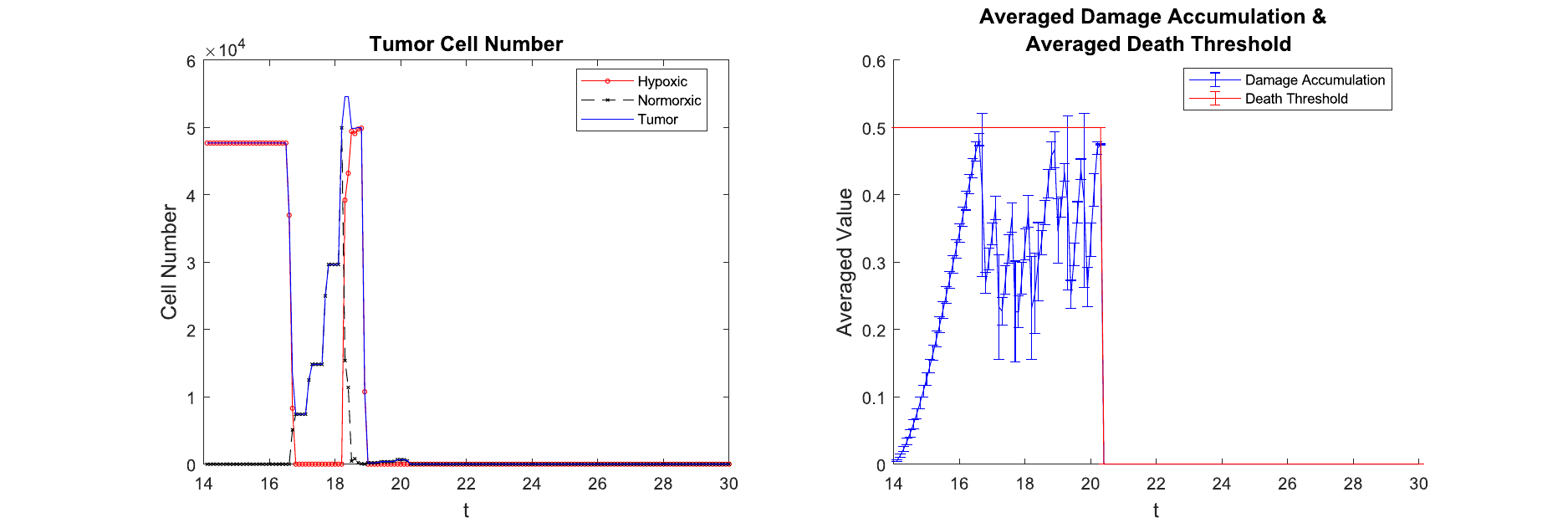}
    \caption{No resistance}
    \label{fig:noresis2}
\end{figure}
The tumour dynamics is illustrated in Figure~\ref{fig:noresis2}.
After the onset of low-dose continuous treatment, the tumour population accumulates cellular damage linearly 
until the damage exceeds the threshold, when the majority of tumour cells die, releasing empty space, and the oxygen level recovers to normoxia. Cell proliferation is then resumed, and after one cell cycle ($t = 0.5$), regrowth occurs and cellular damage increases again. 
During the rebound period, the growth dynamics can be described by the static-increase pattern, with the increase-decrease pattern registering damage accumulation dynamics simultaneously. This can be explained
by the fact that the cell cycles of tumour cells are synchronised and that the newly formed two daughter cells after division inherit half of the damage from their mother cell, thus leading to a sharp decline in average damage. 

The local maximum of the tumour cell population during this regrowth
period is when the average damage approaches the threshold, and the tumour cell population has a sharp decline toward
complete eradication after this peak. A similar oscillation pattern is observed between
hypoxic and normoxic cells during the rebound period.

\subsection{Pre-existing Resistance}
If a small fraction of resistant cells are present before treatment, continuous low-dose treatment is ineffective in removing the tumour, and the tumour population survives.

We define the \textbf{declining point} to be the first time the tumour population begins to decrease, and the \textbf{shifting point} as the first time the tumour composition changes to a completely resistant type, and the death threshold has zero deviations.
The damage clearance rates are set to $p_{r}=0.2, 0.3, 1$, and three growth patterns occur.

\begin{figure}[htbp]
    \centering
    \hspace*{-0.5cm}\includegraphics[width=1.05\linewidth]{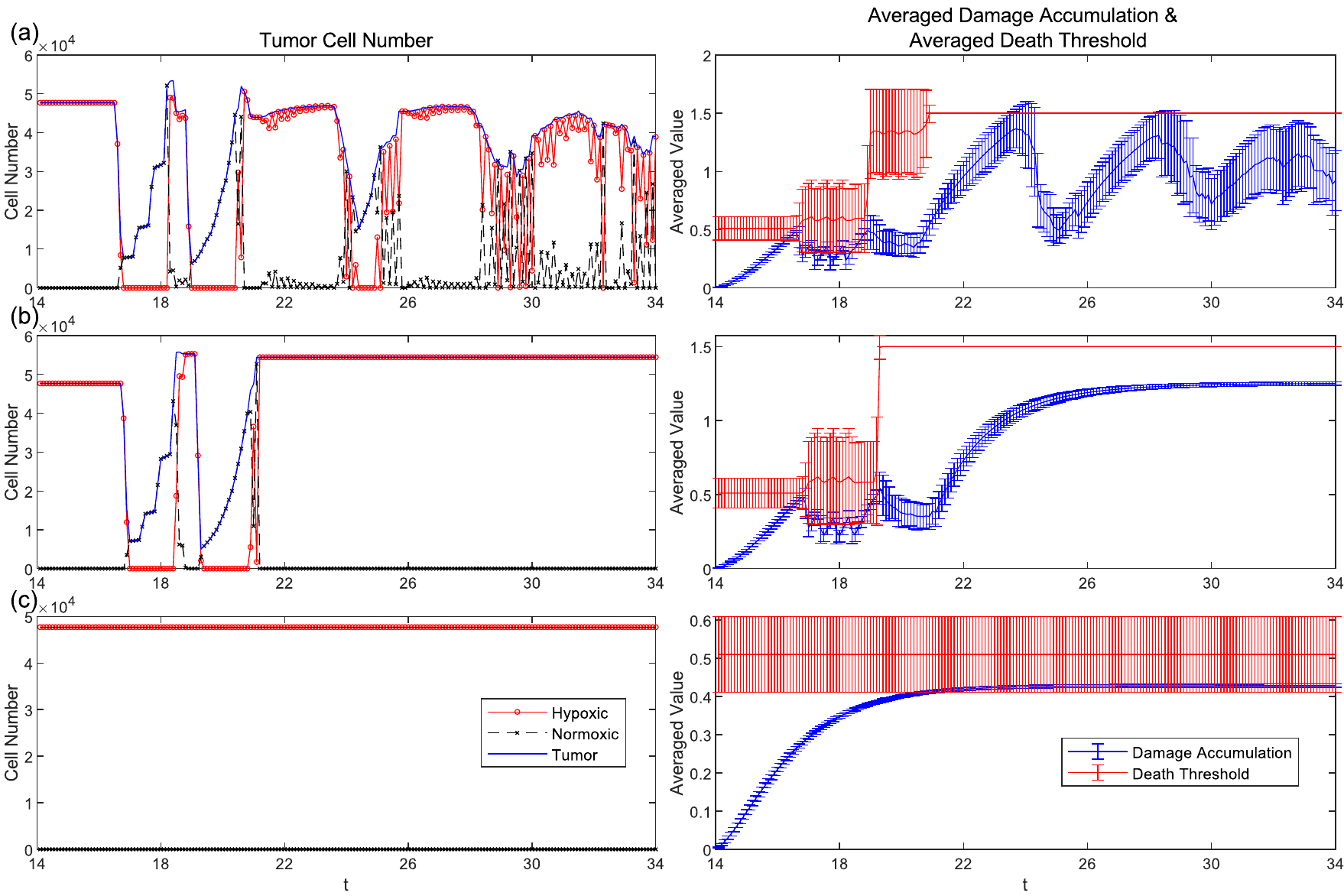}
    \caption{$p_{r} = 0.2, 0.3, 1$}
    \label{fig:preresis1}
\end{figure}

\begin{figure}[htbp]
    \centering
    \hspace*{-2cm}\includegraphics[width=1.05\linewidth]{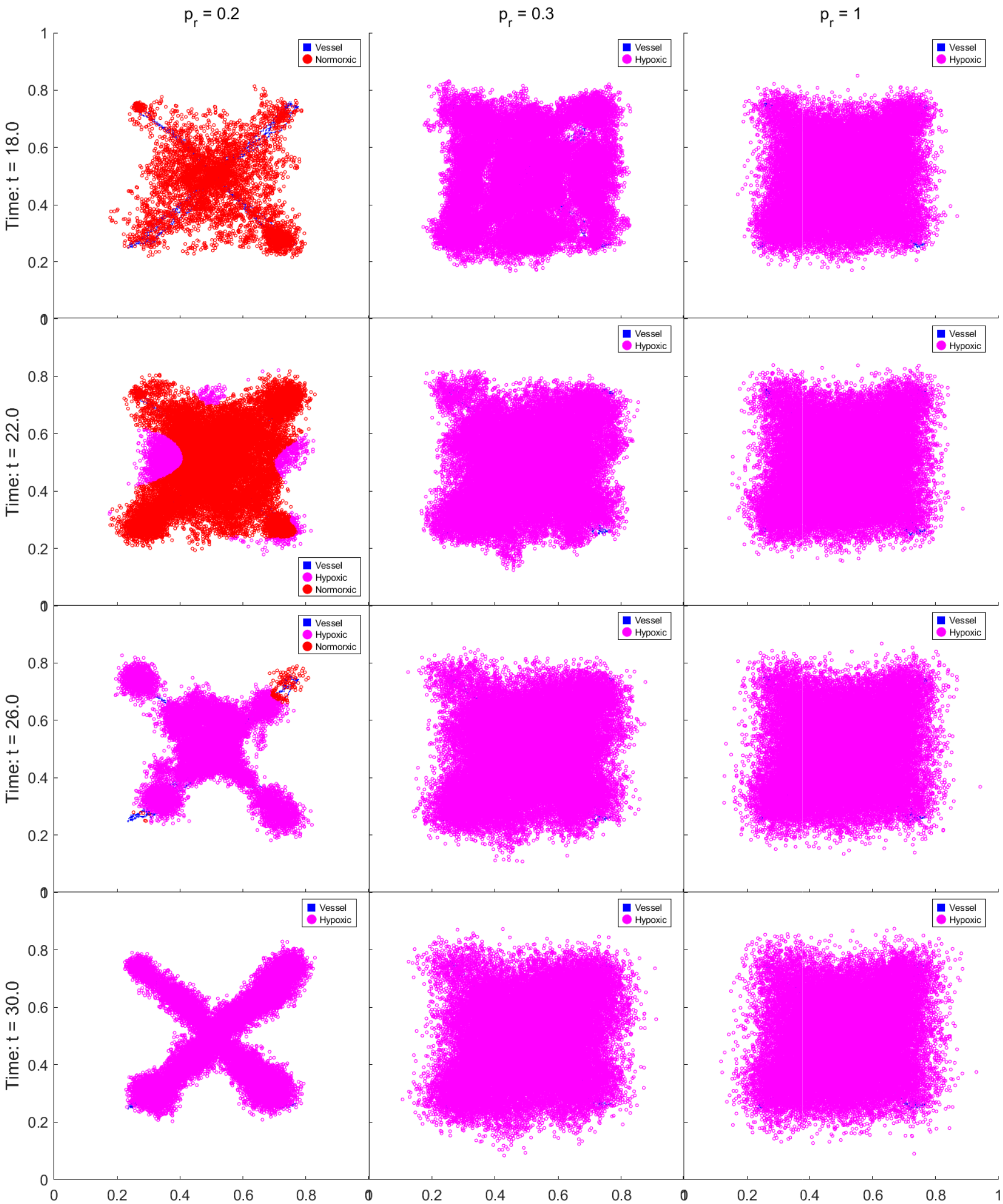}
    \caption{$p_{r} = 0.2, 0.3, 1$}
    \label{fig:preresis2}
\end{figure}

The tumour population is the most vulnerable, with the least damage removed each time step, under $p_{r}=0.2$, and both the tumour population and the average damage feature frequent oscillations after the shifting point. 
The sensitive cells are all killed by the drug, and the average damage is below the threshold, ensuring the survival of resistant cells at a lower level after the shifting point. 

In another extreme case where $p_{r}=1$ with complete damage elimination rendering the treatment ineffective, 
both sensitive and resistant cells survive therapy. The average damage is below the threshold, and ultimately, no fluctuation occurs, and consequently, no tumour apoptosis occurs.

In the intermediate case with a moderate clearance rate $p_{r}=0.3$, there are twice bulk death events followed by tumour replication until they reach the carrying capacity. The surviving clone is resistant. The average damage and tumour population become stable several cycles after the shifting point, with the damage below the threshold.

The mechanisms underlying the protection of the tumour under treatment at low versus high clearance rates are distinct:
In a low $p_{r}$ context, such as $p_{r}=0.2$, persistent oscillations in the damage level are driven primarily by cellular apoptosis, causing the release of vacant space and reducing oxygen consumption. 
Following this, cell division ensues, resulting in symmetric segregation of damage between the two daughter cells, which facilitates a decrease in the damage level.
On the other hand, in scenarios characterised by high clearance rates, for example, when $p_{r}\geq 0.3$, there exists an expedited removal of cellular damage. In contrast to frequent oscillations of damage levels in the low $p_{r}$ case, this rapid clearance leads to a stable increase in average damage levels, whose asymptotic values remain below the threshold. 

Another interesting phenomenon sin Figure~\ref{fig:preresis2} is the increased propensity of tumour cells to aggregate near vessels,
which is the most pronounced in the case of $p_{r}=0.2$. Blood vessels play a dual role in the tumour microenvironment: they supply oxygen to sustain tumour growth while also providing a pathway for drug delivery for tumour eradication.
In this low-dose regimen, the effect of oxygen outweighs the drug, in favour of tumour growth in the adjacent area of the vessels.
In contrast, the drug has a dominating effect elsewhere, leading to the death of tumour cells distant from the vessels.
When the clearance rate increases, the dual role of vessels is not as noticeable as before, and the distribution of tumour cells is more uniform, indicative of the reduced effectiveness of the drug. 

The predominant influence of oxygen around vessels can be ascribed to its ability to maintain cellular proliferation,
which contributes to the halving of the damage level in contrast to a linear progression of drug-induced damage. 
Both the drug and oxygen exert their effects most significantly near vessels: the drug's rapid action is sufficient to eliminate tumour cells within this vicinity; however, concurrently, oxygen facilitates cellular proliferation and ultimately aids in the reduction of damage level between tumour cells. This intricate interplay results in the attenuation of damage, making the region surrounding blood vessels a sanctuary and a survival niche for the tumour population. This result also confirms the prediction of the model in \cite{yang2023multiscale} that the development of angiogenic networks near the tumour promotes a microenvironment that supports tumour growth, thus enhancing drug resistance and increasing the survival rate of tumour cells.

Comparisons between cases exhibiting no resistance with those displaying pre-existing resistance
demonstrates that even a tiny fraction of pre-existing resistant cells can
confer a significant survival advantage in the face of therapeutic interventions. 
The objective of cancer treatment administration is to control the burden of tumours and prolong the lives of patients in the hope of eliminating them and achieving a cure. 
However, if a subset of resistant cells persists, their survival can lead to tumour repopulation with more resistant phenotypes, leading to increased tumour aggressiveness
and reduced treatment response. Therefore, this phenomenon of pre-existing resistance not only increases the likelihood of tumour recurrence but also plays a crucial role in the broader context of treatment failure. 

\subsection{Spontaneous Mutation}
\begin{figure}[htbp]
    \centering
    \hspace*{-2cm}\includegraphics[width=1.25\linewidth]{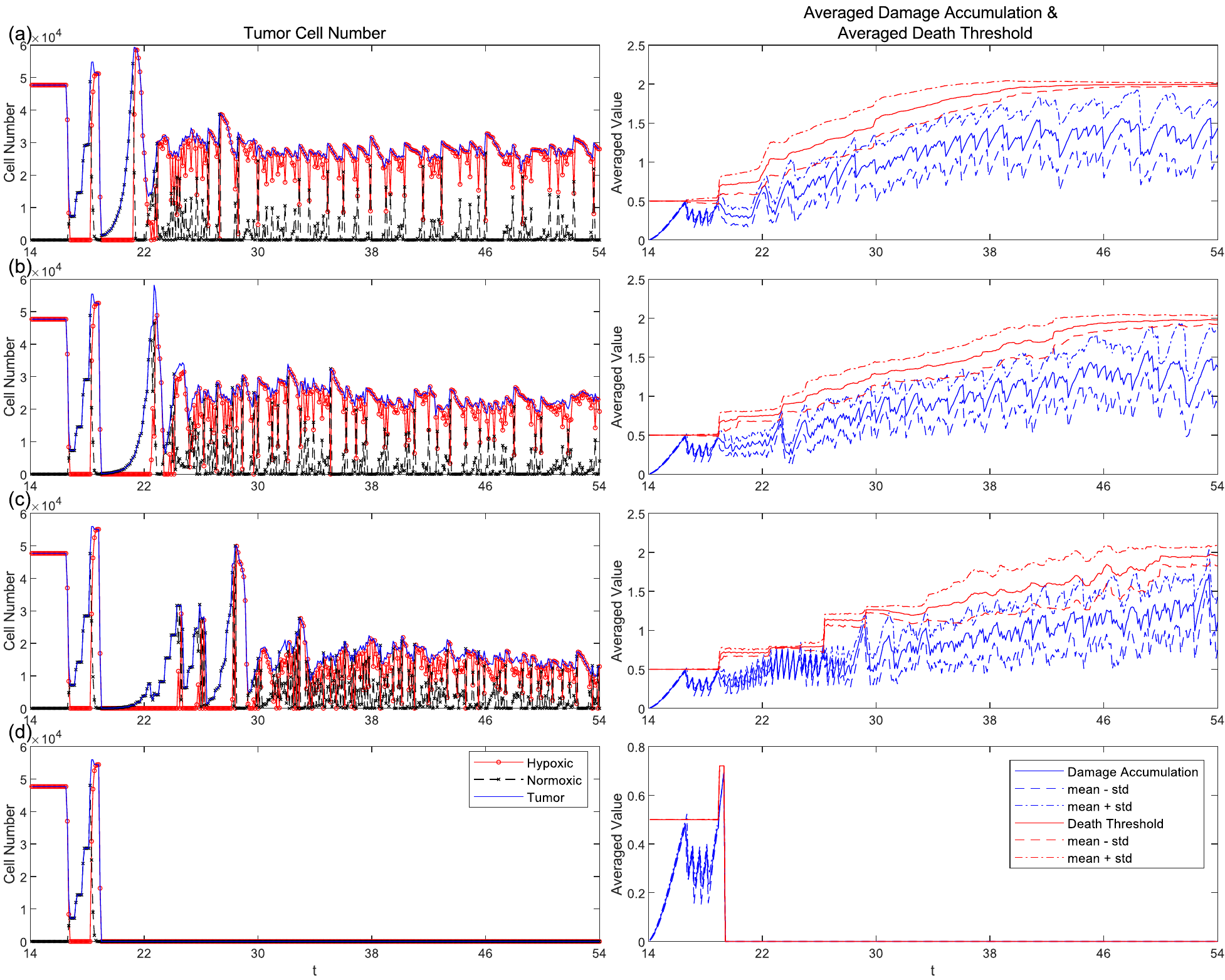}
    \caption{Spontaneous mutation case with $\mu  = 10^{-1}, 10^{-2}, 10^{-3}, 10^{-4}$}
    \label{fig:mutspon1}
\end{figure}

According to \cite{nowak2006evolutionary,araten2005quantitative}, the rate of stem cell mutation per replication
ranges from $10^{-6}-10^{-2}$. We have taken a Poisson distribution for the mutation events, and since the time step is $\Delta t=0.1$, the intensity of the mutation for the distribution takes four values
$\mu =10^{-1}, 10^{-2}$, $10^{-3}$, and $10^{-4}$, while assuming that these mutations occur at the initiation of treatment. 
$p_{r} = 0.2$, so the treatment is successful without mutation events. 

\begin{figure}[htbp]
    \centering
    \hspace*{-2cm}\includegraphics[width=1.25\linewidth]{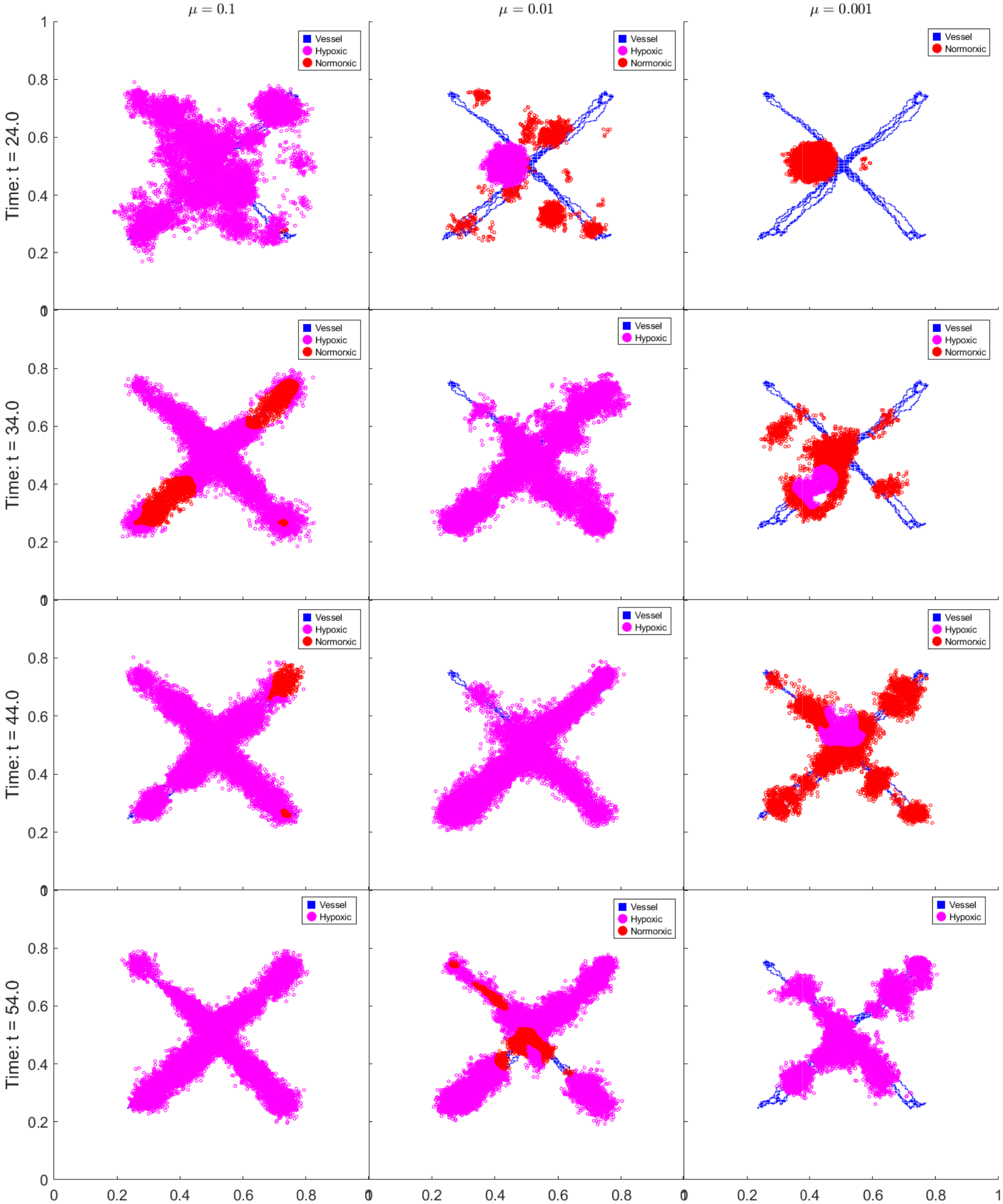}
    \caption{Spantaneous mutation case with $\mu  = 10^{-1}, 10^{-2}, 10^{-3}, 10^{-4}$}
    \label{fig:mutspon3}
\end{figure}

In scenarios where the mutation rate reaches a moderate threshold ($\mu \geq 10^{-3}$), the emergence of treatment progresses slowly, leading to a sustained tumour burden. 
The tumour population experiences a series of events characterised as \textbf{bulk death-rapid increase},
wherein a significant decline in the tumour cell population is followed by a rebound phase marked by a steep increase in cell numbers.
There are two main reasons for these bulk death events:
\begin{enumerate}[(i)]
    \item The first occurs when the average damage inflicted on tumour cells approaches the threshold, accompanied by minimal deviations among individual cells, resulting in a rapid bulk death event.
    \item The second scenario involves considerable variations in damage among tumour cells, such that the upper damage range (i.e., $[\mbox{mean}, \mbox{mean}+\mbox{std}]$) 
overlaps significantly with the lower threshold range (i.e., $[\mbox{mean}-\mbox{std},\mbox{mean}]$). 
\end{enumerate}
With a noticeable portion of the sensitive cells being eliminated, the tumour population composed of more resistant cells
demonstrates an elevated capacity to withstand drug treatment. As a consequence, the tumour population eventually reaches a state of \textbf{ quasi-stabilisation} at a new carrying capacity that is diminished compared to pretreatment levels.
This quasi-stabilisation is characterised by ongoing oscillations in both tumour cell counts and the proportions
of hypoxic and normoxic cells, and the amplitude of such oscillations remains relatively subdued. 
During this phase, there may be a small fraction of cell death, and tumour cells may adapt by raising their death threshold. 
In contrast, with a lower mutation rate ($\mu =10^{-4}$),
the tumour population experiences effective elimination after several cycles of treatment.  
Figure~\ref{fig:mutspon3} provides a visual representation of tumour cells that gather near blood vessels, further supporting previous findings that indicate that the presence of vascular networks can significantly increase tumour cell survival. 

The evolutionary dynamics presented in Figure~\ref{fig:mutspon1} compellingly illustrates how tumour cells
adapt to the selective pressures exerted by drug treatment.
To be precise, the tumour population is initially marked by phenotypic heterogeneity. Drug treatments eliminate phenotypes with a lower survival advantage
while promoting the survival of those that are well-adapted. This selective process diminishes intratumor competition,
which subsequently reduces surrounding nutrient consumption, facilitating the replication of resistant phenotypes.
As the tumour population consists of cells with elevated resistance and improved survival advantages, drug treatment continues to eradicate a minor proportion of the population, forcing the tumour to further elevate its death threshold
and thus diminishing the efficacy of the treatment over time. 

\begin{figure}[htbp]
    \centering
    \hspace*{-2cm}\includegraphics[width=1.25\linewidth]{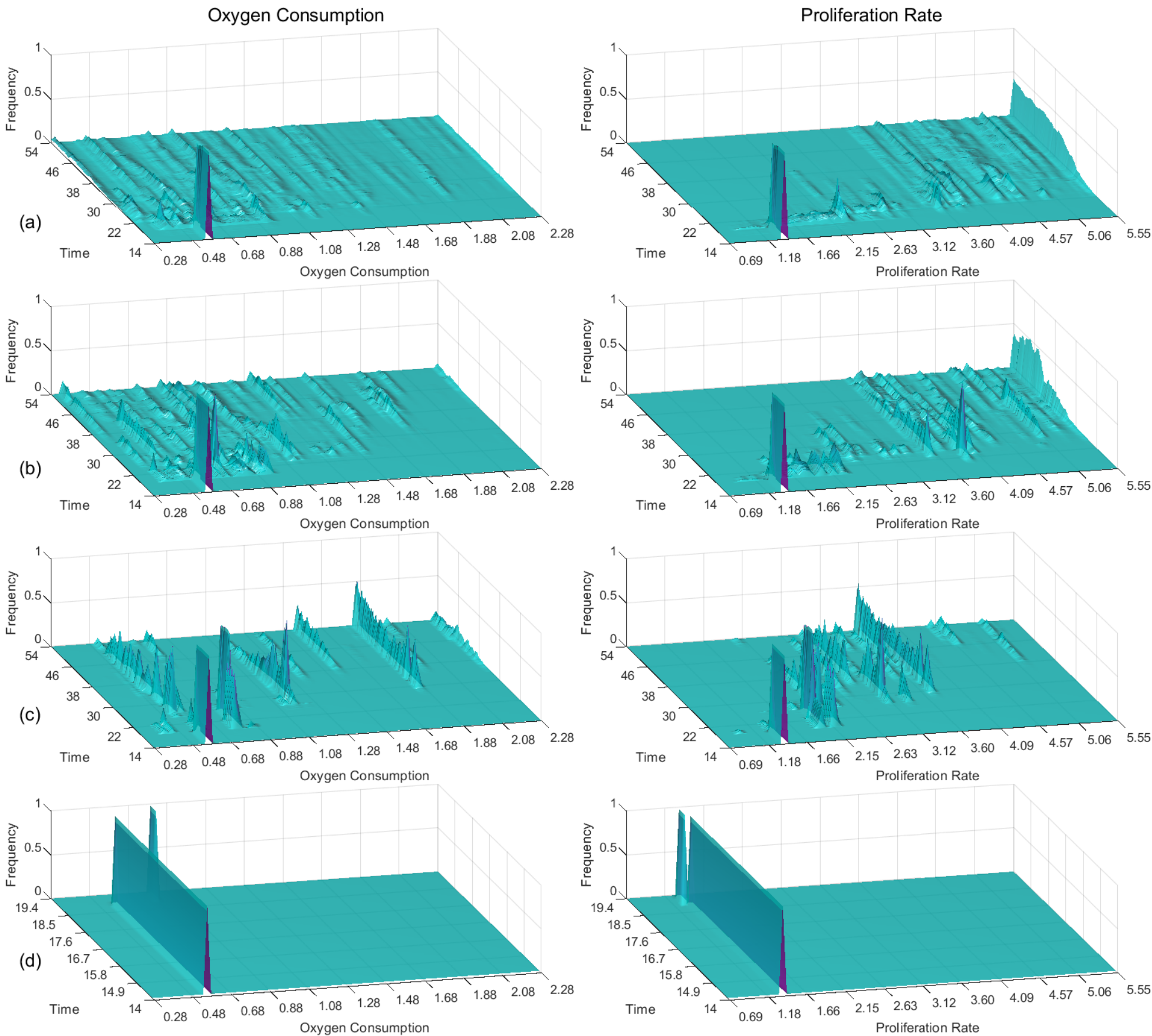}
    \caption{Temporal evolution of oxygen consumption rate and proliferation rate distributions}
    \label{fig:mutspon2}
\end{figure}

Figure~\ref{fig:mutspon2} illustrates the temporal evolution of tumour phenotype distributions in varying mutation rates, using our random mutation algorithm. 
At a very low mutation rate (e.g. $\mu =10^{-4}$), the rarity of mutations results in the distributions of both the oxygen consumption rate and the proliferation rate being tightly concentrated around their initial trait values, $\rho _{o}$ and $\alpha_{n}$, respectively. 
This insufficient mutation frequency does not cause substantial alterations in the composition of the tumour phenotype, thus excluding the development of any protective advantages within the tumour population.
As the mutation rate escalates from moderate ($\mu =10^{-3}$) to elevated levels ($\mu =10^{-1}$), notable shifts are observed: the distribution of oxygen consumption rates transitions from a multimodal
to a more dispersed configuration, compared to the remained unimodal proliferation distribution
whose peak shifts towards phenotypes characterised by expedited proliferation. 
The rigorous environmental constraints imposed by drug treatments and hypoxic conditions
exert considerable selection pressure, ultimately favouring phenotypes that exhibit increased resistance and accelerated proliferation. These findings elucidate the relative importance of diverse tumour phenotypes in enabling adaptation to drug treatment:
\begin{enumerate}[(i)]
\item Given that oxygen concentrations in the simulation do not drop below the apoptosis threshold $o_{apop}$, variations in oxygen consumption rates do not precipitate cell death due to extreme hypoxia. 
\item However, the proliferation rate serves as a critical factor in tumour cell survival. When environmental conditions become conducive to tumour proliferation, characterised by adequate oxygen levels and limited overcrowding, higher proliferation rates at fixed mutation rates
lead to promoted mutation frequencies, thereby boosting the likelihood of producing cells that are
more resistant to treatment. These resistant cells are preferentially selected for survival, 
which subsequently ensures the persistence of rapidly proliferating cells. 
As a result, high proliferation rates and high resistance traits reinforce each other, creating a positive feedback loop. This interplay ultimately leads to the dominance of cells
that possess both high proliferation and high resistance characteristics.
\end{enumerate}
This discussion explains the divergent evolutionary dynamics associated with two distinct traits in tumour cells, which include the oxygen consumption rate and the proliferation rate, emphasising the importance of the proliferation rate as a viral characteristic that allows tumour cells to adapt to pharmacological interventions.
 
A notable observation can be drawn from Figure~\ref{fig:mutspon1} and Figure~\ref{fig:mutspon2}, which illustrate that, between various mutation rates from $10^{-4}$ to $10^{-1}$, the dynamics of the tumour population, the accumulation of damage, and the two frequency distributions exhibit gross qualitative and quantitative similarities
during the initial $10$ cell cycles after the commencement of treatment. This consistency suggests that
the accumulation of mutations requires time to confer any tangible survival advantage to the tumour population. 
In particular, at the end of the tenth cell cycle, specifically from $t=18.9$ to $t=19$, the population numbers decrease to their nadir, with counts changing from $13231$ to $1439$, 
$13765$ to $185$, $16900$ to $49$ and $16390$ to $1$, 
in $\mu =10^{-1}$, $10^{-2}$, $10^{-3}$ and $10^{-4}$, respectively. 
Concurrently, there is a sudden increase in the death threshold and peaks in frequency distributions
pertaining to oxygen consumption rates, and proliferation rates dissipate, exhibiting a sudden transition wherein the proliferation rate tends to shift towards the highest value as the mutation rate increases,
indicating a potential dominance of more resistant and fastest-proliferating cells in survival. 
These abrupt changes observed in the death threshold, the oxygen consumption rate, and the proliferation rate can be understood as follows.
The death threshold is directly related to the survival of the tumour population. 
Only cells that acquire higher resistance through mutations can withstand sudden eradication, leading to the lowest population count at $t=19$. During mutation processes, these resilient cells also develop different rates of oxygen consumption and proliferation than their original trait values. 
The prevalence of these resistant cells, which harbour altered oxygen consumption and proliferation rates, accounts for the sudden loss of peaks in the frequency distributions of both rates.
After $t=19$, the tumour shows a capacity for repopulation in at least moderate mutation rate scenarios, while it may extinguish under the lowest mutation rate condition. Our findings on the qualitatively similar dynamics of tumour cell counts during the early cell cycles highlight the critical time required
for mutations to accumulate and to bestow a survival advantage.
The earlier and more frequently mutations occur, the greater the survival advantage conferred on the tumour population, subsequently leading to a more pronounced decline in treatment efficacy. 
These observations corroborate the conclusions drawn by \cite{yang2023multiscale} regarding the enhanced drug resistance
observed in tumours harboring earlier mutations.

\begin{figure}[htbp]
    \centering
    \hspace*{-2cm}\includegraphics[width=1.25\linewidth]{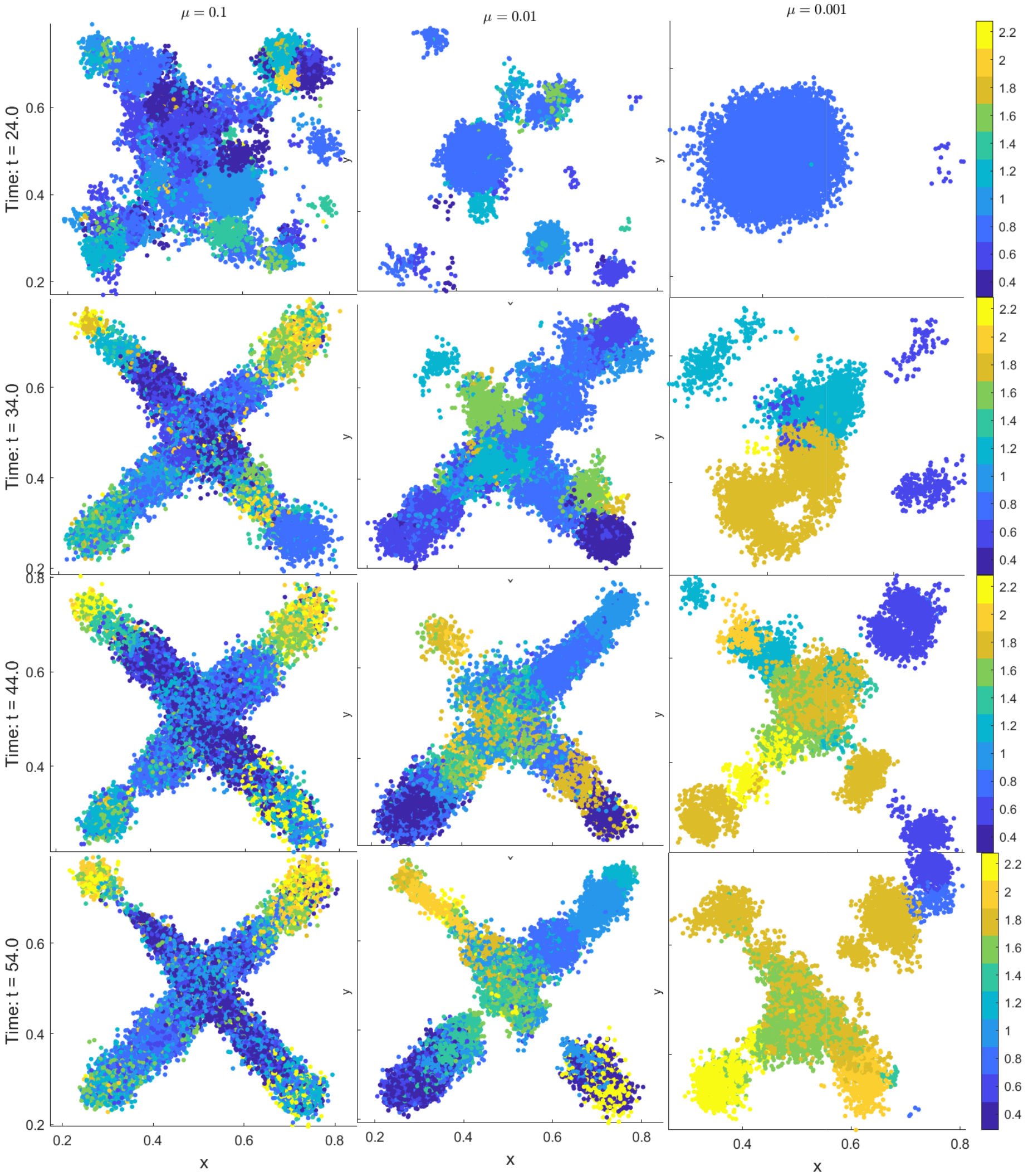}
    \caption{Spatial distribution of oxygen consumption rate}
    \label{fig:mutspon4}
\end{figure}

\begin{figure}[htbp]
    \centering
    \hspace*{-2cm}\includegraphics[width=1.25\linewidth]{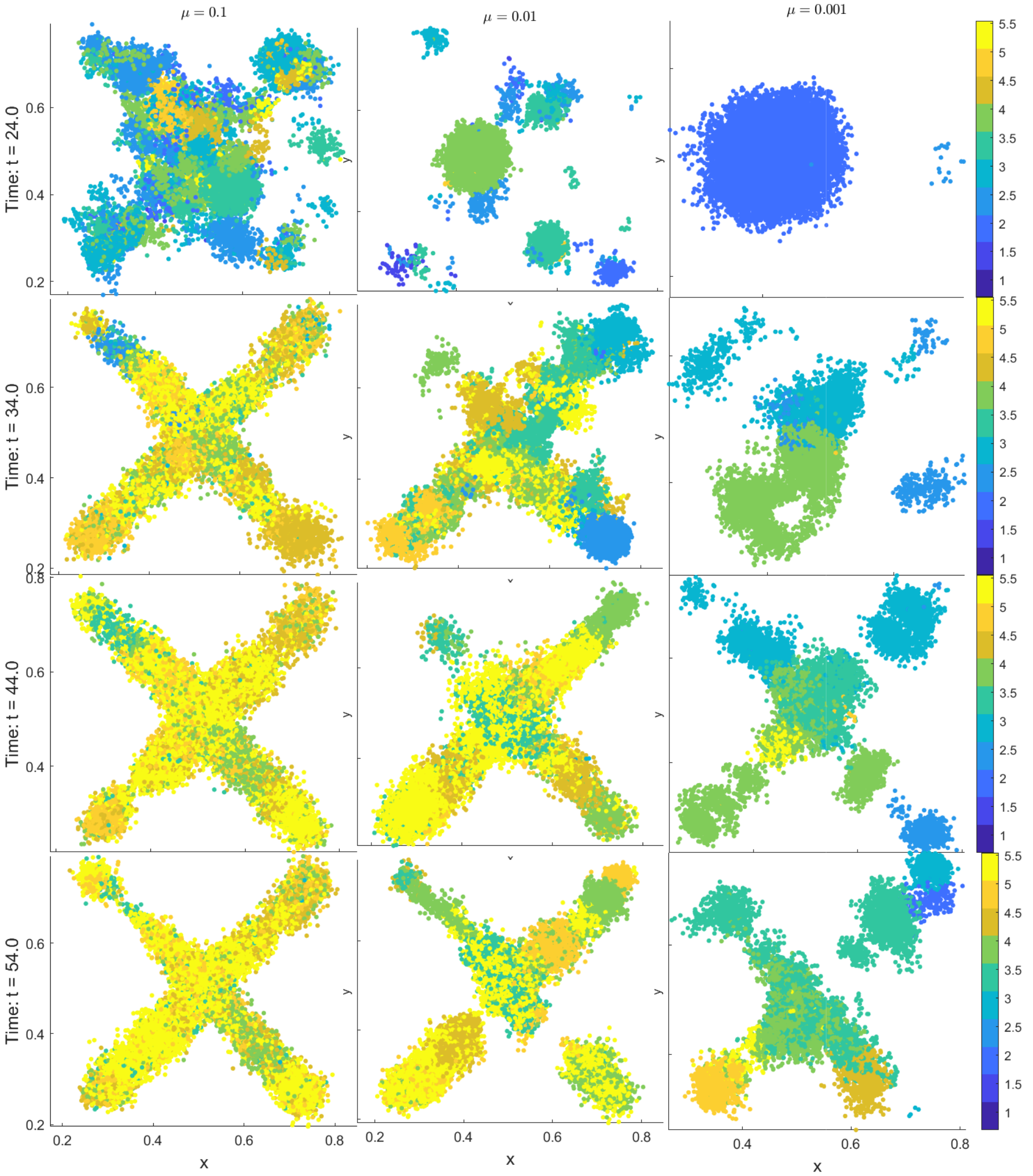}
    \caption{Spatial distribution of proliferation rate}
    \label{fig:mutspon5}
\end{figure}

\begin{figure}[htbp]
    \centering
    \hspace*{-2cm}\includegraphics[width=1.25\linewidth]{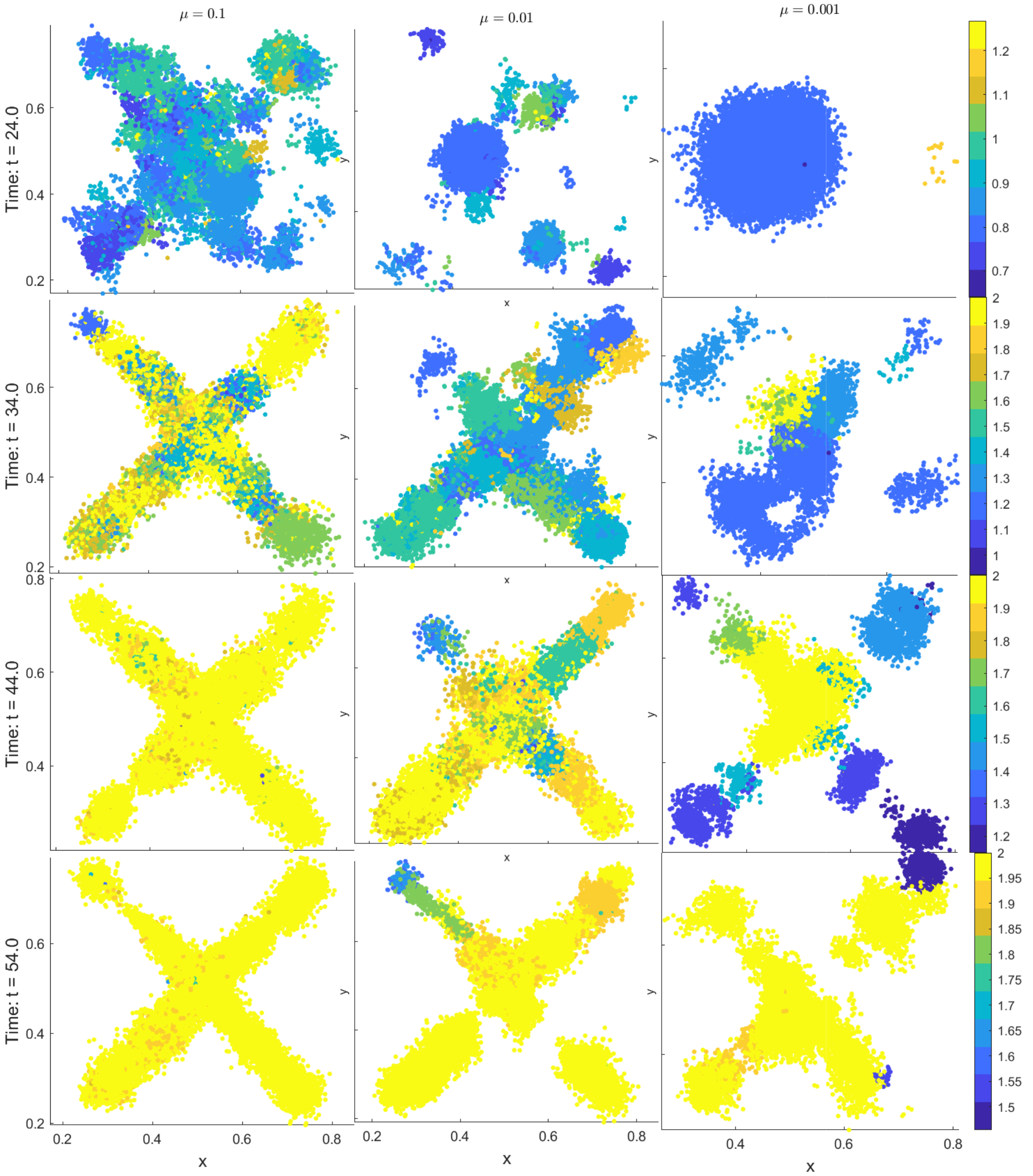}
    \caption{Spatial distribution of resistance trait}
    \label{fig:mutspon7}
\end{figure}

Figure~\ref{fig:mutspon4}, Figure~\ref{fig:mutspon5} and Figure~\ref{fig:mutspon7} depict the spatial displacement of cells
characterised by varying trait values. We do not observe a correlation between oxygen concentration and displacement of
oxygen consumption rates. In contrast, cells that exhibit increased proliferation rates originate primarily from hypoxic regions. 
Furthermore, our findings highlight profound connections between high proliferation rates and enhanced resistance traits: the regions occupied by highly proliferating cells largely overlap with those inhabited by resistant cells. 
In addition, the appearance and prevalence of cells with increased proliferation rates appear to be in sync with those
exhibiting enhanced resistance traits. These connections reaffirm our findings regarding
the mutual reinforcement between high proliferation rates and high resistance traits.

We can classify the resistance mechanisms into two categories:
\begin{enumerate}[(i)]
\item \textbf{Passive adaptation}: In instances of pre-existing resistance, the adverse microenvironment selectively favours resistant cells
through the elimination of sensitive counterparts and the resultant decrease in intratumour competition.
In this paradigm, the trait values of individual tumour cells remain static, yet the overall resistance of the tumour population
shifts due to alterations in its compositional makeup, a process described by passive adaptation.
And survival of the tumour population is guaranteed by the proliferation process, which halves cellular damage.
\item \textbf{Active adaptation}: In contrast, in scenarios involving spontaneous mutations, the capacity to mutate empowers the tumour to
actively adapt to environmental adversities. The detrimental microenvironment imposes selection pressures
that favour resistant and rapidly proliferating phenotypes. In this context, the tumour population increasingly resists treatment, enabling survival even in the presence of severe pharmacological interventions.
This type of adaptation is termed active adaptation, as the tumour population can actively evolve.
\end{enumerate}

These insights have substantial clinical implications for the formulation of effective treatment strategies, as most cancer therapies tend to create a more hostile tumour microenvironment.
For treatment strategies to be effective, it is crucial that the local microenvironment actively contributes to the treatment process, or at least does not hinder it through these selection mechanisms. 
One potential solution is to combine therapies, some of which suppress the selective impact of the microenvironment, while others work to control tumour cell growth cytostatically or eliminate tumour cells cytotoxically. 
Although this promising approach is currently challenging to implement in practice, it is likely to become feasible in the future as the critical role of the microenvironment unravels.

\subsection{Different Treatment Strategies}

As we have mentioned in \ref{sub4.2}, we consider two treatment strategies: continuous and pulsed infusion.
Each treatment period lasts for $t=50$, and we vary the time of treatment and treatment dosage.
We have seven treatment strategies:
\begin{table}[htbp]
     \centering
    \begin{tabular}{|l|l|l|l|l|c|c|c|c|c|}
     \hline
     Treatment& $t_{on}$& $t_{off}$ & $S_{d}$ & strategy & pre-existing & $\mu = 0.1$ & $\mu =0.01$ & $\mu =0.001$ & $\mu =0.0001$    \\
     \hline
     $1$ & $10$ & $40$ & $10$ & pulsed & $\usym{2713}$ & $\usym{2717}$ & $\usym{2717}$ & $\usym{2713}$ & $\usym{2713}$ \\ 
     \hline
     $2$ & $20$ & $30$ & $5$ & pulsed & $\usym{2717}$ & $\usym{2717}$ & $\usym{2717}$ & $\usym{2713}$ & $\usym{2713}$ \\
     \hline
     $3$ & $30$ & $20$ & $10/3$ & pulsed & $\usym{2717}$ & $\usym{2717}$ & $\usym{2717}$ & $\usym{2717}$ & $\usym{2713}$ \\
     \hline
     $4$ & $40$ & $10$ & $5/2$ & pulsed & $\usym{2717}$ & $\usym{2717}$ & $\usym{2717}$ & $\usym{2717}$ & $\usym{2713}$ \\
     \hline
     $5$ & $50$ & $0$ & $2$ & continuous & $\usym{2717}$ & $\usym{2717}$ & $\usym{2717}$ &$\usym{2717}$ & $\usym{2713}$ \\
     \hline
     $6$ & $50$ & $0$ & $5$ & continuous & $\usym{2713}$ & $\usym{2717}$ & $\usym{2713}$ & $\usym{2713}$ & $\usym{2713}$ \\ 
     \hline
     $7$ & $50$ & $0$ & $10$ & continuous & $\usym{2713}$ & $\usym{2713}$ & $\usym{2713}$ & $\usym{2713}$ & $\usym{2713}$ \\
     \hline 
     \end{tabular}
     \caption{Treatment strategies}
     \label{lb2}
\end{table}
The cell damage clearance rate is set to $p_{r}=0.2$, and the mutation rate is set to
$\mu =0.1$.
The results of treatment for pre-existing and spontaneous resistance are shown in Table~\ref{lb2}, where $\usym{2713}$ and $\usym{2717}$
indicate the success and failure of treatment, respectively. 

\section{Biological Implications}
\label{sec:biology}
Non-small cell lung cancer (NSCLC) is the most common type of lung cancer
and is the leading cause of cancer-related deaths in the United States. 
A subset of NSCLC patients that harbour mutations in the epidermal growth factor receptor (EGFR) gene,
including small in-frame deletions in exon $19$ ($19$ dels) and a point mutation within exon $21$ (L$858$R), 
responds well to tyrosine kinase inhibitors (TKI) initially. 
Almost half of patients who develop acquired resistance to TKI are associated with a second site mutation within the exon $20$ of EGFR, termed T$790$M. 

In \cite{godin2007oncogenic,mulloy2007epidermal}, the authors perform surrogate kinase assays in Sf$9$ 
transfectants and transformation assays in fibroblasts, demonstrating that the T$790$M mutation, when combined with the characteristic activating mutations L$858$R or $19$ dels, confers synergistic oncogenic activity. In \cite{chmielecki2011optimization}, paradoxical findings involving the growth disadvantage that the T$790$M mutation confers are raised, based on which the authors anticipate that synergistic oncogenic activity results from genetic alterations other than T$790$M that allow for phenotypic changes to more aggressive ones.

Our findings on the following two aspects in the spontaneous mutation scenario remain in alignment with the above explanation for both the survival advantage and the drug resistance conferred by the T$790$M mutation. 
\begin{enumerate}[(i)]
    \item the earlier and more frequent mutations confer a greater survival advantage upon the tumour population,
    subsequently leading to a more pronounced decline in treatment efficacy;
    \item the mutual reinforcement between high proliferation rates and high resistance traits, including the dominance and the synchronized appearance of cells with 
    high proliferation rates and high resistance traits.
\end{enumerate} 
More frequent mutation events facilitate the emergence of genetic alterations other than T$790$M, and the dominance of highly proliferating and more resistant cells is confirmed by the experimental results in \cite{godin2007oncogenic,mulloy2007epidermal}.

\section{Future Directions}
\label{sec:future}
\begin{enumerate}[(i)]
    \item Evidence has shown that, when applied alone, antiangiogenic factors have not given the expected results. 
    However, when antiangiogenic factors are applied in combination with cytotoxic therapies (chemotherapy and radiation),
    they have proven to reinforce the efficiency of therapies and produce an increase in survival \cite{travasso2011tumor}.
    Many studies on this topic focus only on the mathematical analysis of their proposed models and lack numerical results.
    \item Angiogenesis is also important in the metastasis process, where some tumour cells can escape the primary tumour and enter the bloodstream via newly formed immature and permeable blood vessels to form new tumour masses in distant parts of the body.
    In most cases, the presence of metastases is correlated with tumour malignancy and indicates a poor prognosis for the patient \cite{billy2009pharmacologically}.
    \item The emergence of resistance to EGFR tyrosine kinase inhibitors (TKIs) in non-small cell lung cancer (NSCLC) is multifaceted, with the T$790$M mutation accounting for approximately $50\%$ of cases. However, other mechanisms such as MET amplification and epithelial-to-mesenchymal transition (EMT) also contribute significantly to resistance. 
    MET amplification mechanism is found in $21\%$ of resistant tumors and can occur independently of T$790$M, indicating a distinct pathway of resistance \cite{bean2007met,wang2019met}. 
    EMT can alter cellular characteristics, contributing to resistance and complicating treatment responses \cite{mckenzie2013vivo}.
    This diversity in resistance mechanisms suggests that our current work can be adapted to account for the complexity of tumor evolution by incorporating multiple mechanisms. 
    \item We assume drug release directly from vessel sites at a rate $S_d(t)$ without any restrictions. The limitations of these assumptions are that, in reality, drug diffusion through vessel walls is significantly influenced by the permeability of the vessels themselves. Permeability affects how drugs move from the bloodstream into surrounding tissues, particularly in the context of tumours, where vascular characteristics can vary greatly. Increased permeability of blood vessels enhances drug concentration in tissues, particularly in tumors, as shown in computational simulations \cite{sadipour2023effect}. The permeability of vessel walls is critical for effective drug delivery, as it determines the rate at which drugs can extravasate into the interstitial space \cite{welter2013interstitial}. In vivo studies using models like the developing chicken embryo have demonstrated that manipulating vascular permeability can enhance drug uptake in tumors, as evidenced by increased dextran extravasation in response to specific treatments \cite{pink2012real}. While increased permeability generally facilitates drug delivery, it can also lead to uneven distribution within tumors, complicating treatment efficacy. Understanding these dynamics is essential for optimizing therapeutic strategies.
\end{enumerate} 

\section{Conclusion}
\label{sec:conclusion}
We apply a hybrid discrete-continuous model to study the impact of angiogenesis on the development of spontaneous and induced drug resistance. We use the reaction-diffusion equation to characterise the oxygen and drug field, and they are implemented by the alternating direction implicit method. We employ the agent-based model to simulate tumour and endothelial cells. We find that, in the case of pre-existing resistance, where resistant cells already exist before treatment, the formation of angiogenic networks near the tumour creates a microenvironment that supports tumour survival and enhances drug resistance. 
In the case of spontaneous mutation-induced resistance, where spontaneous mutations lead to drug resistance, 
earlier and more frequent mutations confer a greater survival advantage on the tumour population. 
There is also a mutually reinforcing relationship between a high proliferation rate and high resistance characteristics, including the final dominance of cells with both characteristics in the cell population and the simultaneous emergence of cells with a high proliferation rate and high resistance characteristics. 
This finding explains two conflicting experimental results about a second mutation, T790M, in non-small cell lung cancer (NSCLC): The presence of the T790M mutation conferred resistance and a growth disadvantage, while other experimental results demonstrated resistance and a growth advantage.

Further research can be conducted to explore the combination of antiangiogenic factors with cytotoxic therapies (chemotherapy and radiation) and the reported improvement in the efficiency of therapies.

\newpage
\appendix
\section{Numerical Methods}
\subsection{Discretization of the Continuous Model}\label{app:a1}
Equation (\ref{eq1.1}) can be written as:
\begin{align*}
\dfrac{\partial n}{\partial t} = D_{n}\mathop{{}\bigtriangleup}\nolimits n - \nabla (\chi (c)n)\cdot \nabla c - \chi (c)n \mathop{{}\bigtriangleup}\nolimits c.  
 \end{align*} 
We discretize the continuous equation (\ref{eq1.1}) by the forward Euler finite difference scheme with the same spacing for $x$ and $y$ directions:
\begin{align*}
n_{i,j}^{k+1} = n_{i,j}^{k}P_{0}+n_{i+1,j}^{k}P_{1}+n_{i-1,j}^{k}P_{2}+n_{i,j+1}^{k}P_{3}+n_{i,j-1}^{k}P_{4},
 \end{align*} 
with $[x,y,t]=[i\Delta x,j\Delta x,k\Delta t]$ and $n_{i,j}^{k}=n(i\Delta x,j\Delta x,k\Delta t)$. The coefficients $P_{0},P_{1},P_{2},P_{3},P_{4}$ are proportional to the probabilities of the endothelial cell moving left, right, down, and up, respectively, and they are given by:
\begin{align*}
P_{0} =& 1-\dfrac{4D_{n}\Delta t}{\Delta x^{2}} - \dfrac{\chi (c_{i,j}^{k})\Delta t}{\Delta x^{2}}\left( c_{i+1,j}^{k}+c_{i-1,j}^{k}+c_{i,j+1}^{k}+c_{i,j-1}^{k}-4c_{i,j}^{k} \right), \\ 
P_{1} =& \dfrac{D_{n}\Delta t}{\Delta x^{2}}-\dfrac{\chi (c_{i,j}^{k})\Delta t}{4\Delta x^{2}}\left( c_{i+1,j}^{k}-c_{i-1,j}^{k} \right),    \\ 
P_{2} =& \dfrac{D_{n}\Delta t}{\Delta x^{2}}+\dfrac{\chi (c_{i,j}^{k})\Delta t}{4\Delta x^{2}}\left( c_{i+1,j}^{k}-c_{i-1,j}^{k} \right),    \\   
P_{3} =& \dfrac{D_{n}\Delta t}{\Delta x^{2}}-\dfrac{\chi (c_{i,j}^{k})\Delta t}{4\Delta x^{2}}\left( c_{i,j+1}^{k}-c_{i,j-1}^{k} \right),    \\
P_{4} =& \dfrac{D_{n}\Delta t}{\Delta x^{2}}+\dfrac{\chi (c_{i,j}^{k})\Delta t}{4\Delta x^{2}}\left( c_{i,j+1}^{k}-c_{i,j-1}^{k} \right), 
 \end{align*} 
where $\chi (c) = \dfrac{\chi _{0}}{1+\alpha c} $ is the non-dimensional form of the chemotaxis function. 
 In our treatment, $P_{0}$ is proportional to the probability of remaining stationary, and we take $\Delta t$ small enough so that $P_{0}$ is nonnegative.
 If some of the other four probabilities is negative, say $P_{1}$, and the counterpart $P_{2}\geq 0$, then this means that the probability of moving left is negative, 
 and we set $P_{2} = P_{2} - P_{1}$ and $P_{1} = 0$. On the other hand, if $P_{1}<0$ and $P_{2}<0$, we set $P_{1}=-P_{2}$ and $P_{2}=-P_{1}$.    

 When deciding the movement of the individual endothelial cell, we will normalize $P_{i}$ to ensure that the sum of $P_{i}$ is equal to 1:
\begin{align*}
\tilde{P}_{i} = P_{i}/\left(\sum\limits_{i=0}^{4}P_{j}\right),\quad i=0,1,2,3,4. 
 \end{align*}  

\subsection{ADI method}
For a general nonlinear reaction-diffusion equation, we use the alternating direction implicit (ADI) method to numerically solve it \cite{morton2005numerical}.
We assume the equation has the form \begin{align*}
\dfrac{\partial u}{\partial t} = D \mathop{{}\bigtriangleup}\nolimits u + f(u,x,y,t),
 \end{align*} 
where $D$ is the diffusion coefficient, which could be zero, and $f(u,x,y,t)$ is the reaction term. 
The implicit scheme for the above equation is:
\begin{align*}
\left( 1-\dfrac{D}{2}\dfrac{\Delta t}{\Delta x^{2}}\delta _{x}^{2} \right)U^{k+1/2} &= \left( 1+ \dfrac{D}{2}\dfrac{\Delta t}{\Delta y^{2}}\delta _{y}^{2}\right) U^{k} + \dfrac{\Delta  t}{2}f(U^{k}), \\ 
\left( 1-\dfrac{D}{2}\dfrac{\Delta t}{\Delta y^{2}}\delta _{y}^{2} \right)U^{k+1} &= \left( 1+ \dfrac{D}{2}\dfrac{\Delta t}{\Delta x^{2}}\delta _{x}^{2}\right) U^{k+1/2} + \dfrac{\Delta  t}{2}f(U^{k}),  
 \end{align*} 
where $U^{k+1/2} = u((k+1/2)\Delta t)$ is the intermediate solution at time $(k+1/2)\Delta t$. 
Here, $\delta _{x}^{2}$ and $\delta _{y}^{2}$ are the second-order central difference operators in the $x$ and $y$ directions, respectively.
The ADI method is implicit and is unconditionally stable.

The zero Neumann boundary conditions are dealt with by second-order central differences. 
For example, at the left boundary $x=0$,
we have 
\begin{align*}
\dfrac{U_{1,j}^{k}-U_{-1,j}^{k}}{2\Delta x} = 0,
\end{align*} 
and this gives $U_{-1,j}^{k} = U_{1,j}^{k}$.
The treatment at other boundaries is the same.

\section{Parameter Estimation}
Most parameter values are taken from the literature, and the remaining parameters are estimated by fitting the model to the experimental data. 
We now list the dimensional and non-dimensional parameter values:
\begin{table}[H]
    \centering
    \hspace*{-2cm}\caption{Parameters used in the simulations.}
    \begin{tabular}{|l|l|l|l|l|}
    \hline
    Parameter& meaning &  D-value& ND-Value & Reference\\
    \hline
    $R_{c}$ & cellular radius& $12.5 \mu \mbox{m}$ & 0.005 & \cite{melicow1982three} \\ 
    \hline
    $D_{c}$ & TAF diffusion coeff.& $1.875\times 10^{-1} \dfrac{\mbox{ mm}^{2}}{\mbox{ h}}$ & $0.12$  & \cite{mac2005differential,addison2008simple} \\
    \hline
    $\xi_{c}$ & TAF decay rate & $3.4722\times 10^{-8} \dfrac{1}{\mbox{ s}}$ &  $0.002$ &\cite{billy2009pharmacologically}\\
    \hline
    $\eta $ & TAF production rate & $1.7\times 10^{-22}\dfrac{\mbox{ mol}}{\mbox{ cell}\cdot\mbox{ s}}$& $6.2669\times 10^{3}$ & \cite{hinow2009spatial} \\
    \hline
    $\lambda $ & TAF uptake rate & &  $0.1$  & \cite{anderson1998continuous}\\
    \hline
    $D_{d}$& drug diffusion coeff. &  & $0.5$  &  \\
    \hline
    $\xi_{d}$& drug decay rate & & $0.01$ & \cite{gevertz2015emergence}  \\
    \hline
    $\rho_{d}$& drug uptake rate &  & $0.5$ & \cite{gevertz2015emergence}  \\
    \hline
    $S_{d}$& drug supply rate & & $2$ & \cite{gevertz2015emergence}   \\
    \hline
    $p_{r}$& damage clearance rate &  & $0.2$ & \cite{gevertz2015emergence}   \\
    \hline
    $D_{o}$& oxygen diffusion coeff. & $1\dfrac{\mbox{ mm}^{2}}{\mbox{ h}}$& $0.64$ & \cite{takahashi1966oxygen}  \\ 
    \hline
    $\xi_{o}$ & oxygen decay rate & $4.34\times 10^{-7}\dfrac{1}{\mbox{ s}}$  & $0.025$ & \cite{anderson2007hybrid}  \\
    \hline
    $\rho_{o}$& oxygen uptake rate & $6.25\times 10^{-17}\dfrac{\mbox{ mol}}{\mbox{ cell}\cdot \mbox{ s}}$ & $34.3881$  & \cite{gevertz2015emergence}  \\
    \hline
    $S_{o}$& oxygen supply rate & & $3.5$ & calibrated   \\
    \hline
    $o_{max}$ & maximum oxygen concentration & $6.7\times 10^{-6}\dfrac{\mbox{ mol}}{\mbox{ cm}^{3}}$ & $1$ & \cite{hinow2009spatial} \\ 
    \hline
    $o_{hyp}$& hypoxia threshold  & &$0.25$&\cite{hinow2009spatial}  \\ 
    \hline
    $o_{apop}$& apoptosis threshold & &$0.05$&\cite{hinow2009spatial}  \\ 
    \hline   
    $D_{n}$& endothelial diffusion coeff. & $2\times 10^{-9} \dfrac{\mbox{ cm}^{2}}{\mbox{ s}}$ &  $4.608\times 10^{-4}$ & \cite{anderson1998continuous}\\
    \hline  
    $\chi_{0}$ & chemotaxis coeff.  & $2600 \dfrac{\mbox{ cm}^{2}}{\mbox{ s}\cdot \mbox{ M}}$ &  $0.38$ & \cite{anderson1998continuous}\\ 
    \hline
    $\alpha $ & & &  $0.6$ &\cite{anderson1998continuous}\\ 
    \hline
    $\psi$& threshold branching age &  & $0.5$ & \cite{anderson1998continuous}   \\
    \hline
    $c_{br}$& &   & $1$ & \cite{flandoli2023mathematical}  \\
    \hline   
    $a_{S}^{death}$& sensitive death threshold & & $0.5$ & \cite{gevertz2015emergence} \\ 
    \hline 
    $Th_{multi}$& death threshold ratio&  & $3$ & \cite{gevertz2015emergence}   \\
    \hline
    $\wp _{age}$& cell cycle time&  $U[0.9,1.1] \mbox{ days}$ & $9/16-11/16$  & \cite{picco2024role}  \\
    \hline
    $\alpha _{n}$& proliferation rate &  & $1.0082-1.2323$ &  $\log(2)/\wp _{age}$  \\ 
    \hline
    $F_{max}$& maximum neighbourhood cell number & & $10$ & \cite{flandoli2023mathematical} \\ 
    \hline   
    \end{tabular}
    \label{non-para-estimation}
\end{table} 

   It is difficult to get the values of the dimensional parameters for many parameters, so we consider the non-dimensionalisation of our model.
   We choose a reference diffusion coefficient $D$ such that $\tau =L^{2}/D=16$ h with $L=5 \mbox{ mm}$. 
   The cell cycle time depends on the types of tumour and ranges from $0.25-23.2$ days \cite{chignola1995heterogeneous,demicheli1991exponential,twentyman1980response,sakuma1980cell,wilson1988measurement}, and the authors in \cite{chignola2005estimating} use the Gompertzian model to estimate $0.307-0.596$ days.
   We choose $8$ hours as the cell cycle time, and it is equal to $t=0.5$ on the non-dimensional time scale.  
   The non-dimensionalization for equations (\ref{eq1.1})-(\ref{eq1.2}) is similar to that in \cite{anderson1998continuous}:
   $n, c, t$ is scaled by introducing appropriate reference variables $n_{0}, c_{0}, \tau $: 
   \begin{align*}
   \tilde{n}=\dfrac{n}{n_{0}},\quad \tilde{c}=\dfrac{c}{c_{0}},\quad  \tilde{t}=\dfrac{t}{\tau}.
   \end{align*} 

   The diameters of tumour cells vary depending on the type of tumour under study, but are in the range of $10 - 100 \mu \mbox{m}$ with an appropriate volume of $10^{-9} - 3\times 10^{-8} \mbox{ cm}^{3}$. We take the dimensional value of $n_{0}$ to be $6.4\times 10^{7}\mbox{ cell }\mbox{ cm }^{-3}$ 
   as in \cite{anderson2007hybrid}, and this corresponds to a cell diameter of $25 \mu \mbox{m}$. We therefore take $n_{0}=6.4\times 10^{7}\mbox{ cell}\mbox{ cm}^{-3}$ ($1.6\times 10^{5}\mbox{ cell}\mbox{ cm}^{-2}$ ).    
   $c_{0} = 10^{-10}\mbox{ M}$. 
   By dropping the tildes, we obtain the non-dimensionalised equations:
   \begin{align}
   \dfrac{\partial n}{\partial t}  & = D_{n}\mathop{{}\bigtriangleup}\nolimits n-\nabla \cdot \left( \dfrac{\chi_{0} }{1+\alpha c}n \nabla c  \right),   \\
   \dfrac{\partial c}{\partial t}  & = D_{c}\mathop{{}\bigtriangleup}\nolimits c-\xi _{c}c+\eta \sum\limits_{k} \chi _{C_{k}^{h}}-\lambda nc,
   \end{align} 
   where 
   \begin{align*}
      &\tilde{D}_{n}=\dfrac{D_{n}}{D},\quad \tilde{\chi}_{0}=\dfrac{\chi _{0}c_{0}}{D},\quad \alpha =\dfrac{c_{0}}{k_{1}},\\ 
      &\tilde{D}_{c}=\dfrac{D_{c}}{D},\quad \tilde{\eta }=\dfrac{\eta \tau n_{0}}{c_{0}},\quad \tilde{\lambda }=\lambda \tau n_{0}, \\ 
      &\tilde{\xi}_{c}=\tau \xi_{c},
   \end{align*} 
   
   The non-dimensionalisation for (\ref{eq1.3}) is 
   \begin{align}
   \dfrac{\partial d}{\partial t}=D_{d}\mathop{{}\bigtriangleup}\nolimits d -\xi_{d}d-\rho _{d}d\sum\limits_{k} \chi _{C_{k}}+S_{d}\sum\limits_{j} \chi _{V_{j}},
   \end{align} 
   where, with $d_{0}$ being the reference drug concentration, 
   \begin{align*}
   \tilde{d}=\dfrac{d}{d_{0}},\quad \tilde{D}_{d}=\dfrac{D_{d}}{D},\quad \tilde{\xi}_{d}=\tau \xi_{d},\quad \tilde{\rho }_{d}=\rho _{d}\tau n_{0},\quad \tilde{S}_{d}=\dfrac{S_{d}\tau n_{0}}{d_{0}},
   \end{align*} 
   
   The non-dimensionalisation for (\ref{eq1.4}) is 
   \begin{align*}
   \dfrac{\partial o}{\partial t}=D_{o}\mathop{{}\bigtriangleup}\nolimits o-\xi _{o}o
   -\rho _{o} \sum\limits_{a \in \Lambda _{t}} \chi _{a}(x,t)+S_{o}\left( 1-o \right)\sum\limits_{j} \chi _{V_{j}},  
   \end{align*} 
   where, with $o_{max}$ being the reference oxygen concentration,
   \begin{align*}
   &\tilde{o}=\dfrac{o}{o_{max}},\quad \tilde{D}_{o}=\dfrac{D_{o}}{D},\quad \tilde{\xi }_{o}=\tau \xi _{o}, \\ 
   &\tilde{\rho }_{o}=\dfrac{\rho _{o}\tau n_{0}}{o_{\max}},\quad \tilde{S}_{o}=S_{o}\tau n_{0}.  
   \end{align*} 
   
   In our simulation, the ND values for $D_{o},\eta$, and $\rho_{o}$ are much larger, leading to unreasonable results. 
   Hence, we adjust $\eta = 10^{3}$. We take $\rho_{o}=0.57$, $D_{o}=0.35$ as in \cite{anderson2007hybrid}.

   Our simulation results confirm the substantial influence of the oxygen supply rate $S_{o}$ 
   on the vascular development of tumours. When $S_{o}$ is small, that is, $S_{o}=3$, the tumour will experience hypoxia and cease to expand. In contrast, a high supply rate, that is, $S_{o}=5$, facilitates enough access of tumour cells to oxygen,
causing exponential tumour growth. We choose an intermediate value $S_{o}=3.5$, which triggers the tumour to expand and eventually reach a new carrying capacity. This carrying capacity is determined by $S_{o}$: The surrounding environment with a faster supply rate accommodates more tumour cells, whose proliferation is mainly restricted by the nearby oxygen availability. So, fine-tuning $S_{o}$ in the range $(3, 5)$ leads to the same qualitative behaviour of this vascular development.

  \section{Flowchart of the Algorithm}
  
   \begin{tikzpicture}[node distance=2cm]
   
   \node (pro1) [process] {Update Cell Positions(pro1)};
   \node (proe1) [process, left of =pro1, xshift=-2cm] {Record Tip Cell Trajectories(proe1)};
   \node (proe2) [process, below of =proe1] {Update TAF, Fibronection Concentrations(proe2)};
   \node (proe3) [process, below of =proe2] {Update Angiogenic Network(proe3)};
   \node (proe4) [process, below of =proe3] {Anastomosis(proe4)};
   \node (proe5) [process, below of =proe4] {Update Tip Cell Age(proe5)};
   \node (dece1) [decision, below of =proe5] {Branching Age Threshold?(dece1)};
   \node (dece2) [decision, below of =dece1] {Local Space?(dece2)};
   \node (proe1a) [process, below of =dece2] {Branching(proe1a)};
   \node (proe6) [process, below of =proe1a] {Endothelial Cell Proliferation(proe6)};
   \node (pro2) [process, below of =pro1] {Update Cellular Drug, Oxygen uptake(pro2)};
   \node (pro3) [process, below of =pro2] {Update Drug, Oxygen Concentrations(pro3)};
   \node (dec1) [decision, below of =pro3, right of=pro3, xshift=2cm] {Enough Oxygen Uptake?(dec1)};
   \node (pro1a) [process, below of =pro3] {Increase Age(pro1a)};
   \node (dec4) [decision, below of =dec1] {Low Oxygen?(dec4)};
   \node (pro4a) [process, below of =dec4] {Apoptosis(pro4a)};
   \node (dec2) [decision, below of =pro1a] {Maturation?(dec2)};
   \node (dec3) [decision, below of =dec2] {Local Space?(dec3)};
   \node (pro3a) [process, below of =dec3] {Proliferation(pro3a)};
   \node (pro4) [process, below of =pro3a] {Mutation(pro4)};
   \node (pro5) [process, below of =pro4] {Update Damage(pro5)};
   \node (pro6) [process, below of =pro5] {Update Drug Exposure Time(pro6)};
   \node (dec5) [decision, right of =pro6, xshift=2cm] {Acquire Resistance?(dec5)};
   \node (dec6) [decision, above of =dec5] {Prolonged Drug Exposure(dec6)};
   \node (pro6a) [process, above of =dec6] {Increase Threshold(pro6a)};
   \node (dec7) [decision, above of =pro6a] {Repair Damage; Damage Tolerable?(dec7)};
   \node (new) [startstop, right of =pro1, xshift=2cm] {New Iteration};

   \draw [arrow] (pro1) -- (proe1);
   \draw [arrow] (proe1) -- (proe2);
   \draw [arrow] (proe2) -- (proe3);
   \draw [arrow] (proe3) -- (proe4);
   \draw [arrow] (proe4) -- (proe5);
   \draw [arrow] (proe5) -- (dece1);
   \draw [arrow] (dece1) -- node[anchor=east] {yes} (dece2);
   \draw [arrow] (dece1) -| ++(-2,0) |- ++(0,-6) node[anchor=east, yshift=0.9cm] {no} -| ++(0.3,0) |- (proe6);
   \draw [arrow] (dece2) -- node[anchor=east] {yes} (proe1a);
   \draw [arrow] (dece2) -| ++(-2,0) |- (proe6);
   \draw [arrow] (proe1a) -- (proe6);
   \draw [arrow] (proe6) -| ++(2,0) |- (pro2);
   \draw [arrow] (pro2) -- (pro3);
   \draw [arrow] (pro3) -| (dec1);
   \draw [arrow] (dec1) -- node[anchor=south] {yes} (pro1a);
   \draw [arrow] (dec1) -- node[anchor=east] {no} (dec4);
   \draw [arrow] (dec4) -- node[anchor=east] {yes} (pro4a);
   \draw [arrow] (dec4) -| ++(-2,0) |- ++(0,-7) -| ++(-2,0) |- ++(0,-0.4) -| (pro5);
   \draw [arrow] (pro1a) -- (dec2);
   \draw [arrow] (dec2) -- node[anchor=east] {yes} (dec3);
   \draw [arrow] (dec2) -| ++(2,0) |- ++(0,-7) -| ++(-2,0) |- ++(0,-0.4) -| (pro5);
   \draw [arrow] (dec3) -- node[anchor=east] {yes} (pro3a);
   \draw [arrow] (dec3) -| ++(2,0) |- ++(0,-5) node[anchor=east, yshift=1.9cm] {no} -| ++(-2,0) |- ++(0,-0.4) -| (pro5);
   \draw [arrow] (pro3a) -- (pro4);
   \draw [arrow] (pro4) -- (pro5);
   \draw [arrow] (pro5) -- (pro6);
   \draw [arrow] (pro6) -- (dec5);
   \draw [arrow] (dec5) -- node[anchor=east] {yes} (dec6);
   \draw [arrow] (dec5) -| ++(2,0) |- ++(0,4.9) node[anchor=east, yshift=-1.85cm] {no} -| ++(-2,0) |- ++(0,-0.4) -| (dec7);
   \draw [arrow] (dec6) -- node[anchor=east] {yes} (pro6a);
   \draw [arrow] (dec6) -| ++(2,0) |- ++(0,2.9) -| ++(-2,0) |- ++(0,-0.4) -| (dec7);
   \draw [arrow] (pro6a) -- (dec7);
   \draw [arrow] (dec7) -- node[anchor=east] {no} (pro4a);
   \draw [arrow] (dec7) -| ++(2,0) |- ++(0,12) node[anchor=east, yshift=-4cm] {yes} -| ++(-0.2,0) |- (new);
   \draw [arrow] (new) -- (pro1);
   
   \end{tikzpicture}

\newpage
\bibliographystyle{unsrt}
\bibliography{reference}

\end{document}